\documentclass[10pt]{amsart}
\usepackage{latexsym}
\usepackage{amssymb}
\usepackage{amsmath}
\usepackage{mathrsfs}
\usepackage{bbm}
\usepackage{bbold}
\let\mathbbm\mathbb
\def\AX{\mathop{\fam0 AX}}
\def\CH{\mathop{\fam0 CH}}

\def\Fin{\mathop{\fam0 Fin}\nolimits}
\def\Fin{\mathop{\fam0 Fin}\nolimits}

\def\Fnc{\mathop{\fam0 Fnc}\, }

\def\On{\mathop{\fam0 On}\nolimits}
\def\Ord{\mathop{\fam0 Ord}\, }
\def\Tr{\mathop{\fam0 Tr}\, }

\def\ZF{\mathop{\fam0 ZF}}
\def\ZFC{\mathop{\fam0 ZFC}}

\def\cyc{\mathop{\fam0 cyc}}
\def\dom{\mathop{\fam0 dom}\nolimits}

\def\im{\mathop{\fam0 im}\nolimits}

\def\mix{\mathop{\fam0 mix}}

\def\upwardarrow{\mathord{
  \hbox to 5pt{\hss$\vcenter{\hbox to 2.4pt{\hss$\mathchar"222$\hss}\hrule}\hss$}
}}
\def\downwardarrow{\mathord{
  \hbox to 5pt{\hss$\vcenter{\hrule\hbox to 2.4pt{\hss$\mathchar"223$\hss}}\hss$}
}}

\begin{document}

\title{
WHAT IS BOOLEAN VALUED ANALYSIS?
}

\author{S.~S. Kutateladze}
\address[]{
Sobolev Institute of Mathematics\newline
%\indent 4 Koptuyg Avenue\newline
\indent Novosibirsk, 630090\newline
\indent Russia}
\email{
sskut@member.ams.org
}
\begin{abstract}
This is a brief overview of the basic techniques of Boolean valued analysis.
\end{abstract}
\thanks{An expanded version of a talk at Ta\u\i{}manov's
seminar on geometry, topology, and their applications
in the Sobolev Institute of Mathematics on September 25, 2006.}
%\thanks{Translated from {\it Siberian Electronic Mathematical Reports}, Vol.~3 (2006),
%384--409.}
\thanks{I am grateful to A.~E. Gutman for his numerous subtle improvements
of the preprint of this talk in Russian.}

\maketitle

\begingroup
%\baselineskip=.985\baselineskip

\centerline{\bf1. Introduction}
\setcounter{section}{1}
\bigskip
   The term ``Boolean valued analysis'' appeared within
   the realm of mathematical logic.
    It was Takeuti, a renowned expert in proof theory, who introduced
    the term. Takeuti  defined  Boolean valued analysis in~\cite[p.~1]{Takeuti}
    as ``an application of Scott--Solovay's Boolean valued models
    of set theory to analysis.''
   Vop\v enka  invented similar models at the same time. That is how the question of
   the title receives an answer in zero approximation.
   However, it would be premature to finish at this stage.
   It stands to reason to discuss in more detail  the following three
   questions.

   \subsection{}{\bf Why should we know anything at all about Boolean valued analysis?}
   Curiosity often leads us in science, and oftener we do what we can.
   However we appreciate that which makes us wiser.
   Boolean valued analysis has this value, expanding the limits of our knowledge
   and taking off  blinds from the eyes of the perfect mathematician,
   mathematician {\it  par excellence}. To substantiate this thesis is the main target
   of the exposition to follow.

 \subsection{}{\bf What need the working mathematician know this for?}
   Part of the answer was given above: to become wiser. There is another,
   possibly more important,   circumstance.  Boolean valued analysis
  not only is tied up with many topological and geometrical ideas but also
   provides a technology for expanding the content of the already available
  theorems. Each theorem, proven by the classical means, possesses
  some new  and nonobvious content that relates to ``variable sets.''
  Speaking more strictly, each of the available theorems generates a whole family
  of its next of kin in disguise which is enumerated by all
   complete Boolean algebras or, equivalently,  nonhomeomorphic Stone spaces.

  \subsection{}{\bf What do the Boolean valued models  yield?}
   The essential and technical parts of this survey are devoted to
   answering the question. We will focus on the general methods independent of
   the subtle intrinsic properties of the initial complete Boolean algebra.
   These methods are simple, visual, and easy to apply. Therefore they may be
   useful for the working mathematician. 
   
Dana Scott foresaw the role of Boolean valued models in mathematics and wrote
   as far back as in~1969 (\cite[p.~91]{Scott}): ``We must ask whether there is 
   any interest in these nonstandard models aside from the independence proof; 
   that is, do they have any mathematical interest?
   The answer must be yes, but we cannot yet give a really good argument.''
   Some impressive arguments are available today.

   Futhermore, we must always keep in mind that
   the Boolean valued models were invented in order to simplify
   the exposition of Cohen's forcing~\cite{Cohen}.
   Mathematics is impossible without proof. {\it Nullius in Verba.}
   Therefore, part of exposition is allotted to the scheme of proving
   the consistency of the negation of the continuum hypothesis with the
   axioms of Zermelo--Fraenkel  set theory (with choice) $\ZFC$.
    Cohen was awarded  a~Fields medal in~1966 for this
    final step in settling Hilbert's
    problem No.~1.\footnote{Hilbert~\cite{Problems} considered it plausible that
``{\it as regards equivalence there are, therefore, only two   assemblages of numbers,
  the countable assemblage and the continuum.''}}

\section{\bf Boolean Valued Models}
The stage of Boolean valued analysis
is some Boolean valued model of~$\ZFC$.  To define these models,
we start with a complete
Boolean
algebra.\footnote{A {\it Boolean algebra\/}$B$ is an algebra
with distinct unity $\mathbb 1$ and zero $\mathbb 0$,
over the {\it two-point
field} $2:=\{0,1\}$, whose every element is idempotent.
Given $a,b\in B$, put $a\le b \leftrightarrow ab=a$. The completeness of~$B$
means the existence of the least upper and  greatest lower
bounds of each subset of~$B$.}
For the sake of  comfort we may  view
the elements $\mathbb 0$ and $\mathbb 1$ of the
initial complete Boolean algebra $B$ and the operations on $B$
as an assemblage of some special symbols
for discussing the validity of mathematical propositions.

\subsection{}{\bf A Boolean valued universe.}
Given an ordinal~$\alpha$, put
$$
{\mathbbm V}_{\alpha}^{(B)}:=\{x\mid
\Fnc(x)\ \wedge\ (\exists
\beta)(\beta<\alpha\ \wedge\ \dom (x)
\subset {\mathbbm V}_{\beta}^{(B)})\ \wedge\
\im (x)\subset B\}.
$$
(The relevant information on ordinals is collected below in~7.2.)
In more detail,
$$
\gathered
{\mathbbm V}_0^{(B)}:=\varnothing,
\\
{\mathbbm V}_{\alpha+1}^{(B)}
:=\{x\mid \mbox{$x$ is\ a function,}\ \dom(x)\subset{\mathbbm V}_{\alpha}^{(B)},\
\im(x)\subset~B\},
\\
{\mathbbm V}_{\alpha}^{(B)}:=\bigcup\limits_{\beta<\alpha}{\mathbbm V}_
{\beta}^{(B)}\quad
\mbox{($\alpha$ is\ a limit\ ordinal)}.
\endgathered
$$
The class
$$
{\mathbbm V}^{(B)}:=\bigcup_{\alpha\in \On}{\mathbbm V}_{\alpha}^{(B)}
$$
is a~{\it Boolean valued universe}.
The elements of~${\mathbbm V}^{(B)}$
are~{\it $B$-valued sets}.
Note that
${\mathbbm V}^{(B)}$
consists only of functions. In  particular,
$\varnothing$
is the function with domain~$\varnothing$
and range~$\varnothing$.
Consequently, the  three ``lower'' floors of~${\mathbbm V}^{(B)}$
are composed as follows:
$
{\mathbbm V}_0^{(B)}=\varnothing$, ${\mathbbm V}_1^{(B)}=\{\varnothing\}$, and
${\mathbbm V}_2^{(B)}=\{\varnothing,(\{\varnothing\},b)\mid  b\in B\}.
$
Observe also that
$
\alpha\le \beta \to {\mathbbm V}_{\alpha}^{(B)}\subset {\mathbbm V}_{\beta}
^{(B)}
$
for all ordinals
$\alpha$
and~$\beta$.
Moreover, ~${\mathbbm V}^{(B)}$ enjoys
the
{\it induction principle}
$$
(\forall x\in {\mathbbm V}^{(B)})\,((\forall y\in \dom (x))\ \varphi(y)
\to
\varphi(x))\to (\forall x\in {\mathbbm V}^{(B)})\,\varphi(x),
$$
with
$\varphi$
a~formula of~$\ZFC$.

\subsection{}{\bf Truth values.}
Consider a~formula $\varphi$ of~$\ZFC$, where
 $\varphi=\varphi(u_1,\dots,u_n)$.
If we replace ~$u_1,\dots,u_n$
with~$x_1,\dots,x_n\in {\mathbbm V}^{(B)}$,
then we obtain a particular assertion about
the objects $x_1,\dots,x_n$.
We will try to ascribe to this assertion some
{\it truth value}.
This truth value~$[\![\varphi]\!]$
must be a~member of~$B$.
Moreover, we desire naturally that all theorems of~~$\ZFC$
become ``as valid as possible'' with respect to the new procedure,
so that  the {\it top\/} $\mathbb 1:=\mathbb 1_B$ of $B$, i.e. the {\it unity\/} of~$B$,
serves as the truth value of  a~theorem of~$\ZFC$.

The truth value of a well-formed formula must be determined by
``double'' recursion.
On the one hand, we must induct on the length of
a formula, considering the way it is constructed from the atomic formulas
of the shape
$x\in y$
and $x=y$.
On the other hand, we must define the truth values of
the atomic formulas while $x$ and $y$ range over ${\mathbbm V}^{(B)}$,
on using the recursive construction of the Boolean valued universe.

Clearly, if
$\varphi$
and
$\psi$
are already evaluated formulas of~$\ZFC$, and
$[\![\varphi]\!]\in B$ and
$[\![\psi]\!]\in B$ are the truth values of these formulas
then we must put

$$
\allowdisplaybreaks
\gathered
\,
[\![\varphi\wedge\psi ]\!]:=[\![\varphi{\rm]\!]}\wedge[\![\psi]\!],
\\
[\![\varphi\vee\psi ]\!]:=[\![\varphi{\rm]\!]}\vee[\![\psi ]\!],
\\
[\![\varphi\to \psi ]\!]:=[\![\varphi{\rm]\!]}\rightarrow [\![\psi]\!],
\\
[\![\neg\varphi ]\!]:=\neg [\![\varphi ]\!],
\\
[\![(\forall x)\,\varphi(x)]\!]:=\bigwedge_{x\in {\mathbbm V}^{(B)}}
[\![\varphi(x)]\!],
\\
[\![(\exists x)\,\varphi(x)]\!]:=\bigvee_{x\in {\mathbbm V}^{(B)}}
[\![\varphi(x)]\!],
\endgathered
$$
where the right-hand sides contain the Boolean operations that correspond to
the logical connectives and quantifiers on the left-hand sides:
$\wedge$
is the meet,
$\vee$
is the join,
$\neg $
is the complementation,
and the implication
$\rightarrow$
is defined by the rule
$a\rightarrow b:=\neg a\vee b$
for $a,b\in B$.
Only these definitions ensure the value ``unity''
for the truth values of the classical tautologies.

We turn now to evaluating the atomic formulas
$x\in y$
and
$x=y$
for
$x,y\in {\mathbbm V}^{(B)}$.
The intuitive idea behind the definition
states that every
$B$-valued set~$y$ is ``fuzzy,'' i.e.
``it contains each member
$z$  of~$\dom (y)$ with probability~$y(z)$.''
Bearing this in mind and attempting to preserve
the logical tautology
$x\in y\leftrightarrow (\exists z\in y)\,(z=x)$
alongside the axiom of extensionality, we arrive to
the following recursive definition:
$$
\gathered
\,
[\![x\in y]\!]:=\bigvee_{z\in \dom(y)} y(z)\wedge[\![z=x]\!],
\\
[\![x=y]\!]:=\biggl(\bigwedge_{z\in \dom (x)} x(z)
\rightarrow[\![z\in y]\!]\biggr)\wedge
\biggl(\bigwedge_{z\in \dom (y)} y(z)\rightarrow[\![z\in x]\!]\biggr).
\endgathered
$$
\noindent
We can now ascribe  some meaning to the formal
records like
$\varphi(x_1,\dots,x_n)$,
where
$x_1,\dots,x_n\in {\mathbbm V}^{(B)}$
and
$\varphi$
is a formula of~$\ZFC$. In other words, we are able to
define  the strict sense in which
the set-theoretic expression
$\varphi(u_1,\dots,u_n)$
is valid for
$x_1,\dots,x_n\in {\mathbbm V}^{(B)}$.

Namely, we will say that a~{\it formula
$\varphi(x_1,\dots,x_n)$
is satisfied inside~${\mathbbm V}^{(B)}$\/}
or  the collection of {\it ~$x_1,\dots,x_n$
possesses the property~$\varphi$
inside~${\mathbbm V}^{(B)}$\/}
provided that
$[\![\varphi(x_1,\dots,x_n)]\!]=\mathbb 1$.
In this event we write
${\mathbbm V}^{(B)}\models \varphi(x_1,\dots,x_n)$.
If some formula $\varphi$ of~$\ZFC$  is expressed in
the natural language of discourse then
by way of pedantry we use quotes:
 ${\mathbbm V}^{(B)}\models\,$``$\varphi$.''
 The {\it satisfaction mark\/} $\models$ implies
 the usage of the model-theoretic expressions of the kind
 ``${\mathbbm V}^{(B)}$~is a~ Boolean valued model
 for~$\varphi$''  instead of~${\mathbbm V}^{(B)}\models\varphi$.

All axioms of the first-order predicate calculus
are obviously valid inside~${\mathbbm V}^{(B)}$.
In particular,

{\bf(1)}
$[\![x=x]\!]=\mathbb 1$,

{\bf(2)}
$[\![x=y]\!]=[\![y=x]\!]$,

{\bf(3)}
$[\![x=y]\!]\wedge[\![y=z]\!]\le [\![x=z]\!]$,

{\bf(4)}
$[\![x=y]\!]\wedge[\![z\in x]\!]\le [\![z\in y]\!]$,

{\bf(5)}
$ [\![x=y]\!]\wedge[\![x\in z]\!]\le [\![y\in z]\!]$.

\noindent
Observe that
$
{\mathbbm V}^{(B)}\models
x=y\wedge\varphi(x)\to\varphi(y)
$
for every formula $\varphi$, i.e.

{\bf(6)}
$[\![x=y]\!]\wedge[\![\varphi(x)]\!]\le [\![\varphi(y)]\!]$.

\section{\bf Principles of Boolean Valued Analysis}

The equality
$[\![x=y]\!]=\mathbb 1$
holding in the Boolean valued universe
${\mathbbm V}^{(B)}$
implies in no way that
the functions
$x$
and~$y$ coincide
in their capacities of members of~${\mathbbm V}$.
For instance, the function vanishing at one of the levels
~${\mathbbm V}_{\alpha}^{(B)}$,
with
$\alpha\ge 1$,
plays the role of the empty set inside~${\mathbbm V}^{(B)}$.
This circumstance is annoying in a few constructions we use in the sequel.

\subsection{}{\bf The separated  universe.}
To overcome this nuisance, we will pass from~${\mathbbm V}^{(B)}$
to the {\it separated Boolean valued universe\/}
$\overline{{\mathbbm V}}^{(B)}$,
keeping for the latter the previous symbol~${\mathbbm V}^{(B)}$;
i.e. we put
${\mathbbm V}^{(B)}:=\overline{{\mathbbm V}}^{(B)}$.
Furthermore, to define
$\overline{{\mathbbm V}}^{(B)}$,
we consider the equivalence
$\{(x,y)\mid [\![x=y]\!]=\mathbb 1 \}$
on~${\mathbbm V}^{(B)}$.
Choosing a member (representative of the least rank)
in each class of equivalent functions,
we come in fact to the separated universe~$\overline{{\mathbbm V}}^{(B)}$.

Note that the implication
$$[\![x=y]\!]=\mathbb 1\to [\![\varphi(x)]\!]=[\![\varphi(y)]\!]$$
holds for every formula $\varphi$~of~$\ZFC$ and members~$x$,~$y$
of~${\mathbbm V}^{(B)}$.
Therefore, we may evaluate formulas over the separated universe,
disregarding the choice of representatives.
Dealing with the separated universe in the sequel,
we will consider, for the sake of convenience, a particular representative
of a class of equivalence rather than the whole class itself
(in contrast to the conventional practice of analysis in treating
the function spaces of the Riesz scale).

\subsection{}{\bf Properties of a Boolean valued model.}
The main properties of a Boolean valued universe~${\mathbbm V}^{(B)}$
are collected in the following three propositions:

{\bf(1)}{\scshape~Transfer Principle.}
 All theorems of~$\ZFC$ are valid inside~${\mathbbm V}^{(B)}$;
i.e., in symbols,
$$
{\mathbbm V}^{(B)}\models \ZFC.
$$

The transfer principle is established by
a bulky check that the truth value of every axiom of~$\ZFC$
is~$\mathbb 1$ and all inference rules increase truth values.
The transfer principle reads sometimes as follows:
``${\mathbbm V}^{(B)}$ is a~Boolean valued model of~$\ZFC$.''
That is how the term ``Boolean valued model of set theory''
enters the realms of mathematics.

{\bf(2)}{\scshape~Maximum Principle.}
To each formula~$\varphi$ of~$\ZFC$
there is a~member~$x_0$ of~${\mathbbm V}^{(B)}$ satisfying
$
[\![(\exists x)\,\varphi(x)]\!]=[\![\varphi(x_0)]\!].
$

In particular, it is true inside~${\mathbbm V}^{(B)}$
that if $\varphi(x)$ for some~$x$ then there
exists a~member~$x_0$ of~${\mathbbm V}^{(B)}$
(in the sense of~${\mathbbm V}$)
satisfying
$[\![\varphi(x_0)]\!]=\mathbb 1$.
In symbols,
$${\mathbbm V}^{(B)}\models (\exists x)\,\varphi(x)\to (\exists x_0)\,
{\mathbbm V}^{(B)}\models \varphi(x_0).$$

\noindent
The maximum principle  means that
$(\exists x_0\in {\mathbbm V}^{(B)})\,[\![\varphi(x_0)]\!]=
\bigvee_{x\in {\mathbbm V}^{(B)}} [\![\varphi(x)]\!]$
for every formula~$\varphi$ of~$\ZFC$. The last formula illuminates
the background behind the term ``maximum principle.''
The proof of the principle consists
in an easy application of the following fact:

{\bf(3)}{\scshape~Mixing Principle.}
Let
$(b_{\xi})_{\xi\in \Xi}$
be a~{\it partition of~unity\/}
in~$B$; i.e.,
a family of members of~$B$
such that
$$
\bigvee_{\xi\in \Xi} b_{\xi}=\mathbb 1,\quad
(\forall \xi,\eta\in \Xi)\,(\xi\ne \eta\to b_{\xi}\wedge b_{\eta}=\mathbb 0).
$$

For every family~$(x_{\xi})_{\xi\in \Xi}$
of elements of~~${\mathbbm V}^{(B)}$
and every partition of unity~$(b_{\xi})_{\xi\in \Xi}$
there is a unique mixing~$(x_{\xi})$
by~$(b_{\xi})$ (or with probabilities $(b_{\xi})$); i.e.,
a~member~$x$
of the separated universe~${\mathbbm V}^{(B)}$
satisfying
$b_{\xi}\le [\![x=x_{\xi}]\!]$
for all~$\xi\in \Xi$.
The {\it mixing\/}
$x$
of a family~$(x_{\xi})$
by~$(b_{\xi})$
is denoted as follows:
$
x=\mix_{\xi\in \Xi}(b_{\xi}x_{\xi})=\mix \{
b_{\xi}x_{\xi}\mid \xi\in \Xi \}.
$
%\noindent
Mixing is connected with the main
particularity of a~Boolean valued model,
the procedure  of  collecting the highest
truth value ``stepwise.''

{\bf(4)}
Take $x\in{\mathbbm V} ^{(B)}$ and $b\in B$.  Define the function
$bx$ on~$\dom(bx)\!:=\dom(x)$ by the rule: $bx:t\mapsto b\wedge x(t)$ for $t\in\dom(x)$.
Then $bx\in{\mathbbm V} ^{(B)}$. Moreover, for all $x,y\in{\mathbbm V}^{(B)}$
we have
 $ [\![x\in by]\!]=b\wedge [\![x\in y]\!]$ and
 $[\![bx=by]\!]=b\rightarrow [\![x=y]\!]$.

 The checking of the first equality
 consists in the straightforward calculation
 of truth values and application of the infinite distributive law:
   $$
 \gathered\ [\![x\in
 by]\!]=\bigvee\limits_{t\in\dom(by)} (by)(t)\wedge [\![t=x]\!]=
 %\\ =
 b\wedge\bigvee\limits_{t\in\dom(y)} y(t)\wedge
 [\![t=x]\!]=b\wedge [\![x\in y]\!].
 \endgathered
 $$
 The second identity uses the rules of inference
 of~2.2:
 $$
 \gathered\ [\![bx=by]\!]=
 %\\ =
\biggl(\bigwedge\limits_{t\in\dom(by)} (by)(t)\rightarrow [\![t\in
 bx]\!]\biggr)\wedge\biggl(\bigwedge\limits_{t\in\dom(bx)} (bx)(t)\rightarrow
[\![t\in by]\!]\biggr)
 \\[0.5\jot]
 =\biggl(\bigwedge\limits_{t\in\dom(y)}\! \left(b\wedge
 y(t)\right)\!\rightarrow\! \left(b\wedge [\![t\in x]\!]\right)\biggr)\wedge
 %\\ \wedge
\biggl( \bigwedge\limits_{t\in\dom(x)}\! \left(b\wedge
 x(t)\right)\!\rightarrow\! \left(b\wedge [\![t\in y]\!]\right)\biggr)
 \\[0.5\jot]
 =\bigwedge\limits_{t\in\dom(y)} \left((b\wedge y(t))\rightarrow
 b\right)\wedge \left((b\wedge y(t))\rightarrow [\![t\in x]\!]\right)
 \\[0.5\jot]
 \wedge\bigwedge\limits_{t\in\dom(x)} \left((b\wedge
 x(t))\rightarrow b\right)\wedge \left((b\wedge x(t))\rightarrow [\![t\in y]\!]\right)
 \\[0.5\jot]
 =\biggl(\bigwedge\limits_{t\in\dom(y)} b\rightarrow (y(t)\rightarrow
 [\![t\in x]\!])\biggr)\wedge \biggl(\bigwedge\limits_{t\in\dom(x)} b\rightarrow (x(t)\rightarrow [\![t\in y]\!])\biggr)
\\ =
 b\rightarrow [\![x=y]\!].
 \endgathered
 $$

 {\bf(5)}
  Given  $b\in B$ and $x\in{\mathbbm V} ^{(B)}$, observe that
  $
 [\![bx=x]\!]=b\vee[\![x=\varnothing ]\!]$ and $[\![bx=\varnothing
 ]\!]=\neg b\vee[\![x=\varnothing ]\!] $.  In~particular,
 ${\mathbbm V} ^{(B)}\models bx=\mix\{bx,\neg b\varnothing\}$.

{\bf(6)}
Let $(b_\xi)$ be a~partition of unity in~$B$ and let
a family $(x_\xi)\subset{\mathbbm V} ^{(B)}$ be such that
 ${\mathbbm V}^{(B)}\models x_\xi\ne x_\eta$
 for all $\xi\ne\eta$. Then there is a~member $x$ of~${\mathbbm V} ^{(B)}$
 satisfying
 $[\![x=x_\xi ]\!]=b_\xi $ for all $\xi$.
Indeed, put $x\!:=\mix(b_\xi x_\xi)$ and
 $a_\xi\!:=[\![x=x_\xi ]\!]$. By hypothesis
 $ a_\xi\wedge a_\eta =[\![x=x_\xi ]\!]\wedge [\![x_\eta
 =x]\!]\le\neg[\![x_\xi\ne x_\eta ]\!]={\mathbb 0}
 $
 for $\xi\ne\eta$. Moreover, $b_\xi\le a_\xi$ for all $\xi$ by
 the properties of mixing. Hence, $(a_\xi)$ is
 a~partition of unity in~$B$ too. On the other hand,
 $
 \neg b_\xi =\bigvee\nolimits_{\eta\ne\xi}
 b_\eta\le\,\bigvee\nolimits_{\eta\ne\xi}a_\eta =\neg a_\xi
 $
 and so $\neg b_\xi\le \neg a_\xi$, i.~e.  $b_\xi\ge a_\xi$. So,
 $(b_\xi)$ and $(a_\xi)$ present the same partition of unity.

 {\bf(7)}
Consider some formulas $\varphi(x)$ and $\psi (x)$ of~$\ZFC$.
Assume that
 $[\![\varphi(u_0)]\!]={\mathbb 1}$ for some $u_0\in{\mathbbm V} ^{(B)}$.
 Then
  $$
   \gathered\
 [\![(\forall\,x)(\varphi(x)\to\psi
 (x))]\!]=\bigwedge \big\{[\![\psi (u)]\!]\mid u\in{\mathbbm V}^{(B)},\,[\![\varphi(u)]\!]={\mathbb 1}\big\},
 \\[1.5\jot]
 [\![(\exists\,x)(\varphi(x)\wedge\psi (x))]\!]=
 \bigvee\big\{[\![\psi (u)]\!]\mid u\in{\mathbbm V}^{(B)},\,[\![\varphi(u)]\!]={\mathbb 1}\big\}.
 \endgathered
 $$
\noindent
These formulas reveal the  ``mosaic'' mechanisms of truth verification in~$\mathbbm V^{(B)}$.

\subsection{}{\bf Functional realization.}
Let $Q$ be the Stone space of a~complete Boolean algebra~$B$.
Denote by $\mathfrak U$  the (separated) Boolean valued
universe~${\mathbbm V}^{(B)}$.
Given $q\in Q$, define the equivalence $\sim_q$
on the class $\mathfrak{U}$ as follows:
 $
 u\sim_q v\ \leftrightarrow\ q\in[\![u=v]\!].
 $
Consider the bundle
 $$
 V^Q:=\big\{\big(q,{\sim_q}(u)\big) \mid q\in Q,\, u\in\mathfrak{U}\big\}
 $$
 and agree to denote the pair
 $\big(q,{\sim}_q(u)\big)$ by $\widehat u(q)$. Clearly,
 for every $u\in\mathfrak{U}$ the mapping
 $\widehat u:q\mapsto \widehat u(q)$ is a section of~$V^Q$.
 Note that to each $x\in V^Q$ there are $u\in\mathfrak{U}$
 and $q\in Q$ satisfying $\widehat u(q)=x$.
Moreover, we have $\widehat u(q)=\widehat v(q)$
 if and only if $q\in[\![u=v]\!]$.

Make each fiber~$V^q$ of~$V^Q$ into an algebraic system of
 signature~$\{\in\}$ by letting
 $
 V^q\models x\,{\in}\, y \ \leftrightarrow\  q\in[\![u\,{\in}\, v]\!],
 $
where $u,v\in\mathfrak{U}$ are such that $\widehat u(q)=x$ and $\widehat v(q)=y$.
Clearly, this definition is sound.
Indeed, if $\widehat u_1(q)=  x$ and $\widehat v_1(q)=y$
 for another pair of elements $u_1$ and~$v_1$ then
 the claims of the memberships $q\in[\![u\,{\in}\, v]\!]$ and $q\in[\![u_1\,{\in}\,
 v_1]\!]$ are equivalent.

It is easy to see that the class of the sets
 $\{\widehat u(A)\mid u\in\mathfrak U\}$,
 with $A$ a~clopen subset of~$Q$, is a base for some topology
 on~$V^Q$. This enables us to view~$V^Q$
 as a~continuous bundle called  a~{\it continuous polyverse}.
By a~{\it continuous section\/} of~$V^Q$  we mean a section
that is a continuous function.  Denote by $\mathfrak C$
 the class of all continuous sections of~$V^Q$.

 The mapping $u\mapsto\widehat u$ is a bijection between~$\mathfrak U$
 and $\mathfrak C$ in which we can clearly find grounds
 for a convenient functional realization of the
 Boolean valued universe ${\mathbbm V}^{(B)}$.
 The details of this universal construction due to
 Gutman and Losenkov  are presented in~\cite[Ch.~6]{IBA}.

\subsection{}{\bf Socialization.}
Mathematics of the twentieth century exhibited many examples
of the achievements obtained by socialization of objects and problems;
i.e., their inclusion into a class of similar objects
or problems.\footnote{Hilbert said in his report~\cite{Problems}:
``In dealing with mathematical problems, specialization plays, as I
believe, a still more important part than generalization. Perhaps in
most cases where we seek in vain the answer to a question, the cause
of the failure lies in the fact that problems simpler and easier than
the one in hand have been either not at all or incompletely solved.
All depends, then, on finding out these easier problems, and on
solving them by means of devices as perfect as possible and of
concepts capable of generalization. This rule is one of the most
important levers for overcoming mathematical difficulties and it seems
to me that it is used almost always, though perhaps unconsciously.''}
Boolean valued models gain a natural status
within category theory.  The idea of a variable set
becomes a cornerstone of the categorical analysis of logic
which is accomplished in topos theory~\cite{Johnstone}.

\subsection{}{\bf Distant modeling.}
 The properties of a~Boolean valued universe reflects
 a~new conception of modeling which may be referred to as
 {\it extramural\/} or {\it distant  modeling}. We will
 explain the essence of this conception by
 comparing it with the traditional approaches.

Encountering two classical models of
the same theory, we usually  seek for a bijection between the
universes of the models.
If we manage to establish such a bijection so as
the predicates and operations translate faithfully
from one model to the other then we speak
about isomorphism between the models. Consequently, this conception of isomorphism
implies an explicit comparison of the models which consists in
witnessing some bijection between their universes of discourse.

Imagine that we are physically unable to compare the models
pointwise. Happily, we take an opportunity
to exchange information with the owner of the other model
by using some means of communication, e.g., by having long-distance
calls.
While communicating, we easily learn that our
partner uses his  model
to operate on some objects that are the namesakes of ours,
i.e., sets, membership, etc.
Since we are interested in $\ZFC$, we ask him
whether or not the axioms of $\ZFC$ are satisfied in his
model.  Manipulating the model,
he returns a positive answer.
After checking that he uses the same inference rules as we do,
we cannot help but acknowledge his model to
be a model of the theory we are all investigating.
It is worth noting that this conclusion still leaves unknown
for us the  objects that make up his universe and
the procedures he uses to distinguish between
true and false propositions about these objects.

 The novelty of distant modeling
 resides in our refusal to identify the universes
 of objects and admittance of some a priori unknown procedures of
 verification of propositions.

\subsection{}{\bf Technology.}
 To prove the relative consistency of
 some set-theoretic propositions we use a~
 Boolean valued universe ${\mathbbm V}^{(B)}$ as follows:
 Let
 $\mathscr T$ and $\mathscr S$ be some enrichments of
 Zermelo--Fraenkel theory $\ZF$ (without choice). Assume that
 the consistency of~$\ZF$ implies the consistency of~$\mathscr S$.
 Assume further that we can define~$B$ so that $\mathscr S\models$
 ``$B$ is a ~complete Boolean  algebra'' and $\mathscr S\models[\![\varphi
 ]\!]={\mathbb 1}$ for every axiom~$\varphi$ of~$\mathscr T$.
 Then the consistency of~$\ZF$ implies the consistency of~$\mathscr T$.
 That is how we use ${\mathbbm V}^{(B)}$ in foundational studies.

 Other possibilities for applying~${\mathbbm V}^{(B)}$
 base on the fact that irrespective of the  choice of a~Boolean algebra~$B$,
 the universe is an arena for testing
 an arbitrary mathematical event.
 By the principles of transfer
 and maximum, every ${\mathbbm V}^{(B)}$  has the objects that play the
 roles of numbers, groups, Banach spaces, manifolds, and whatever
constructs of mathematics that are already introduced into practice or
still remain undiscovered.  These objects may be viewed as some nonstandard
 realizations of the relevant originals.

 All celebrated and not so popular theorems acquire interpretations
 for the members of~${\mathbbm V}^{(B)}$, attaining the top
 truth value. We thus obtain a~new technology of comparison
 between  the interpretations of mathematical facts in the universes over
 various  complete Boolean algebras.  Developing the relevant tools
 is the crux of Boolean valued analysis.

\section{\bf Ascending and Descending}
No comparison is feasible without some dialog between
~${\mathbbm V}$ and~${\mathbbm V}^{(B)}$.
We need some sufficiently convenient
mathematical toolkit for the comparative analysis of the
interpretations of the concepts and facts of mathematics
in various models.
The relevant {\it technique of ascending and descending\/} bases
on the operations of the canonical embedding, descent, and ascent
to be addressed right away.
We start with the canonical embedding of the von Neumann
 universe.

\subsection{}{\bf The canonical embedding.}
Given $x\in {\mathbbm V}$,
denote
 by~$x^{\scriptscriptstyle\wedge}$
the {\it standard name\/} of~$x$
in~${\mathbbm V}^{(B)}$;
i.e., the member of~${\mathbbm V}^{(B)}$  that is defined by
recursion as follows:
$$
\varnothing^{\scriptscriptstyle\wedge}:=\varnothing,\quad \dom (x^{\scriptscriptstyle\wedge}):=
\{y^{\scriptscriptstyle\wedge}\mid y\in x \},
\quad \im (x^{\scriptscriptstyle\wedge}):=\{\mathbb 1\}.
$$
Note a few properties of the mapping
$x\mapsto x^{\scriptscriptstyle\wedge}$
which we will use in the sequel.

{\bf(1)}
If
$x\in {\mathbbm V}$
and $\varphi$ is a~formula of~$\ZFC$ then
$$
\gathered
\,
[\![(\exists y\in
x^{\scriptscriptstyle\wedge})\,\varphi(y)]\!]={%\mboxstyle
\bigvee}\{[\![\varphi(z^{\scriptscriptstyle\wedge})]\!]\ |\
z\in x \},
\\
[\![(\forall y\in x^{\scriptscriptstyle\wedge})\,\varphi(y)]\!]
={%\mboxstyle
\bigwedge}\{[\![\varphi(z^{\scriptscriptstyle\wedge})]\!]\ |\
z\in x \}.
\endgathered
$$

{\bf(2)}
If
$x$ and $y$ are members of~${\mathbbm V}$
then we see by transfinite induction that
$$
%\gather
x\in y\leftrightarrow {\mathbbm V}^{(B)}\models x^{\scriptscriptstyle\wedge}\in y^{\scriptscriptstyle\wedge},
%\\
%$$
%$$
\quad
x= y\leftrightarrow {\mathbbm V}^{(B)}\models x^{\scriptscriptstyle\wedge}=y^{\scriptscriptstyle\wedge}.
%\endgather
$$
In other words, the standard name may be viewed as
the embedding functor from~${\mathbbm V}$ to~${\mathbbm V}^{(B)}$.
Put ${\mathbb 2}\!:=\{{\mathbb 0},{\mathbb 1}\}$, where
$\mathbb 0$ and $\mathbb 1$ are  the zero and unity of~$B$.
This two-point algebra is acknowledged as one of the
hypostases of~2. Clearly,
${\mathbbm V}^{(\mathbb 2)}$ serves naturally as a~subclass of~${\mathbbm V}^{(B)}$.
Undoubtedly, the standard name sends
${\mathbbm V}$ onto
${\mathbbm V}^{(2)}$,
which is distinguished by the following proposition:

{\bf(3)}
$
(\forall u\in {\mathbbm V}^{(\mathbb 2)})\,(\exists !\, x\in {\mathbbm V})\
{\mathbbm V}^{(B)}\models u=x^{\scriptscriptstyle\wedge}.
$

{\bf(4)}
A formula is {\it restricted\/}
provided that each bound variable enters
this formula under the sign of some restricted
quantifier; i.e., a quantifier
ranging over a~set.
Strictly speaking, each bound variable must be restricted by
a quantifier of the form~$(\forall x\in y)$
or~$(\exists x\in y)$
for some~$y$.

{\scshape Restricted Transfer Principle.}
{\sl Given a restricted formula $\varphi$ of~$\ZFC$ and
a~collection $x_1,\dots,x_n\in {\mathbbm V}$,
the following holds:}
$$
\varphi(x_1,\dots,x_n)\leftrightarrow {\mathbbm V}^{(B)}\models
\varphi(x_1^{\scriptscriptstyle\wedge},\dots,x_n^{\scriptscriptstyle\wedge}).
$$
Working further in the separated  universe~$\overline{{\mathbbm V}}^{(B)}$,
we agree to preserve the symbol~$x^{\scriptscriptstyle\wedge}$
for the distinguished  member of the class corresponding to~$x^{\scriptscriptstyle\wedge}$.

{\bf(5)}
By way of example, we mention the
useful consequence of the restricted transfer principle:
$$
\gathered
\mbox{``$\Phi$ is a~correspondence from $x$ to $y$''}
%\leftrightarrow
\\
\leftrightarrow
{\mathbbm V}^{(B)}\models
\mbox{``$\Phi^{\scriptscriptstyle\wedge}$
is a correspondence from $x^{\scriptscriptstyle\wedge}$
to
$y^{\scriptscriptstyle\wedge}$''};
\\
%$$
%$$
\mbox{``$f$ is a function from $x$ to $y$''
$\leftrightarrow $
${\mathbbm V}^{(B)}\models$
``$f^{\scriptscriptstyle\wedge}$ is a function from $x^{\scriptscriptstyle\wedge}$
to
$y^{\scriptscriptstyle\wedge}$''}
\endgathered
$$
(in this event,
$f(a)^{\scriptscriptstyle\wedge}=f^{\scriptscriptstyle\wedge}(a^{\scriptscriptstyle\wedge})$
for all $a\in x$).
Therefore, the standard name may be considered as a
covariant functor from the category of sets
and mappings (or correspondences) inside~${\mathbbm V}$
to an appropriate subcategory~${\mathbbm V}^{(\mathbb 2)}$
inside the separated universe~${\mathbbm V}^{(B)}$.

{\bf(6)}~A set~$x$
is  {\it finite\/} provided that
$x$
coincides with the range of a function
on a~finite ordinal. We express this symbolically as
$\Fin(x)$. Hence,
$
\Fin(x):=(\exists\,n)(\exists\,f)(n\in\omega\wedge
\Fnc(f)\wedge\dom(f)=n\wedge\im(f)=x)
$
where  $\omega:=\{0,1,2,\dots\}$ as usual.
The above formula  clearly fails to be restricted.
However, the transformation of finite sets proceeds happily
under the canonical embedding.
Indeed, let
$\mathscr P_{\Fin}(x)$
stand for the class of finite subsets of~$x$; i.e.,
$
\mathscr P_{\Fin}(x):=\{y\in\mathscr P(x)\mid\Fin (y)\}.
$
Then
$
{\mathbbm V} ^{(B)}\models\mathscr P _{\Fin} (x)^{\scriptscriptstyle\wedge} =
\mathscr P _{\Fin} (x^{\scriptscriptstyle\wedge})$
for all $X$.
Since the formula specifying the powerset of~a~set~$x$ is unrestricted, in general we may claim only that
$[\![\mathscr P(x)^{\scriptscriptstyle\wedge}\subset
\mathscr P(x^{\scriptscriptstyle\wedge})]\!]=\mathbb 1$.

{\bf(7)}
 Let $\rho$ be an~automorphism of a~complete  Boolean algebra~$B$,
 and let $\psi_\rho$ be the~member of~${\mathbbm V}^{(B)}$
 defined by the function
 $\{(b^{\scriptscriptstyle\wedge},\rho (b))\mid b\in B\}$.
 Then

 \begin{itemize}

 \item[\bf(a)]
 $\rho (b)=[\![b^{\scriptscriptstyle\wedge}\in\psi_\rho ]\!]$
 for all $b\in B$;

 \item[\bf(b)]
  $[\![A^{\scriptscriptstyle\wedge}\subset\psi_\rho
 \to\left(\bigwedge A\right)^{\scriptscriptstyle\wedge}\in\psi_\rho ]\!]={\mathbb 1}$
  for $A\subset B$ if and only if $\rho\left(\bigwedge A\right)=\bigwedge\rho (A);$

 \item[\bf(c)]
 $[\![\psi_\rho$ is an~ultrafilter on~$B^{\scriptscriptstyle\wedge}]\!] ={\mathbb 1}$.

\end{itemize}

{\bf(8)}
Let $\pi$ be a~homomorphism of~$B$ to a complete Boolean algebra~$C$.
By recursion on the relation $y\in\dom(x)$
 we can define the mapping
 $\pi^*:{\mathbbm V} ^{(B)}\to{\mathbbm V} ^{(C)}$
 such that $\dom(\pi^* x)\!:=\{\pi^*y \mid y\in\dom(x)\}$ and
 $
 \pi^*x:v\mapsto\bigvee\big\{\pi(x(z))\mid z\in\dom(x),\,\pi^*z=v\big\}.
 $
If $\pi$ is injective then so is~$\pi ^*$. Moreover,
$\pi ^*x:\pi ^*y\mapsto\pi (x(y))$ for all $y\in\dom(x)$.

Let $\mathfrak U$ be an~ultrafilter on a~Boolean algebra $B$, and
 let $\mathfrak U\,'$ be the dual ideal; namely,
  $\mathfrak U\,'\!:=\{\neg b\mid \,b\in\mathfrak U\}$. Then the factor-algebra
   $B/\mathfrak U\,'$ has two members and we may identify it with
 the two-point  Boolean algebra ${\mathbb 2}\!:=\{{\mathbb 0},{\mathbb 1}\}$.
 The factor-homomorphism $\pi:B\to{\mathbb 2}$ is not complete in general;
 i.e., $\pi$ fails to preserve all joins and meets.
 This prevents us from interconnecting the
 truth values on~${\mathbbm V}^{(B)}$ and ${\mathbbm V} ^{({\mathbb 2})}$. However, if~$\pi$ is
 complete (i.e.,  $\mathfrak U$ is a principal ultrafilter), then for
  every formula $\varphi(x_1,\dots,x_n)$ and every tuple
  $u_1,\dots,u_n\in{\mathbbm V} ^{(B)}$ we have
 $$
 {\mathbbm V} ^{({\mathbb 2})}\models
 \varphi(\pi ^*u_1,\dots,\pi ^*u_n)\leftrightarrow
 [\![\varphi(u_1,\dots,u_n)]\!]\in\mathfrak U\, ,
 $$
 since  the formulas $\pi (b)={\mathbb 1}$ and
 $b\in\mathfrak U$ are equivalent for~$b\in B$.

{\bf(9)}
 Given an ultrafilter $\mathfrak U$, we can use factorization to construct
 a~model that differs from~${\mathbbm V} ^{({\mathbb 2})}$.
 To this end, we introduce in~${\mathbbm V}^{(B)}$ the relation
  ${\sim _{\mathfrak U}}$ by the rule
 $
 {\sim_{\mathfrak U}}:=\{(x,y)\in{\mathbbm V} ^{(B)}\times{\mathbbm V}^{(B)}\mid
 [\![x=y]\!]\in\mathfrak U\}.
 $
Clearly,  ${\sim _{\mathfrak U}}$ is an equivalence on~${\mathbbm V}^{(B)}$.
Denote by ${\mathbbm V} ^{(B)}/\mathfrak U$
the factor-class of ${\mathbbm V} ^{(B)}$ by
 ${\sim _{\mathfrak U}}$ and furnish it with the relation
  $
 {\in_{\mathfrak U}}:=\{(\widetilde x,\widetilde y)\mid x,y\in{\mathbbm V}^{(B)},\,
 [\![x\in y]\!]\in\mathfrak U\},
 $
 where $x\mapsto\widetilde x$ is the canonical factor-mapping
 from ${\mathbbm V}^{(B)}$ to~${\mathbbm V}^{(B)}/\mathfrak U$.
 We can show that
  $
 {\mathbbm V}^{(B)}/\mathfrak U\models
 \varphi(\widetilde x_1,\dots,\widetilde x_n)\leftrightarrow
 [\![\varphi(x_1,\dots,x_n)]\!]\in\mathfrak U
 $
 for $x_1,\dots,x_n\in{\mathbbm V}^{(B)}$ and a formula $\varphi$
 of~$\ZFC$.

\subsection{}{\bf Descent.}
Given an arbitrary member~$x$
of a~(separated) Boolean valued universe~${\mathbbm V}^{(B)}$,
define the {\it descent\/}
$x{\downarrow}$ of~$x$ as follows:
$$
x{\downarrow}:=\{y\in {\mathbbm V}^{(B)}\mid [\![y\in x]\!]=\mathbb 1 \}.
$$
We list the simplest properties of descending:

{\bf(1)}
The class~$x{\downarrow}$
is a~set; i.e.,
$x{\downarrow}\in {\mathbbm V}$
for all~$x\in {\mathbbm V}^{(B)}$.
If
$[\![x\ne \varnothing]\!]=\mathbb 1$
then
$x{\downarrow}$
is a~nonempty set.

{\bf(2)}
Let
$z\in {\mathbbm V}^{(B)}$
and
$[\![z\ne \varnothing]\!]=\mathbb 1$.
Then for every formula~$\varphi$
of~$\ZFC$ we have
$$
\gathered
\,
[\![(\forall x\in z)\,\varphi(x)]\!]
=\bigwedge\{[\![\varphi(x)]\!]\mid x\in z{\downarrow} \},
\\
[\![(\exists x\in z)\,\varphi(x)]\!]
=\bigvee\{[\![\varphi(x)]\!]\mid x\in z{\downarrow} \}.
\endgathered
$$
Moreover, there is some
$x_0\in z{\downarrow}$
satisfying
$[\![\varphi(x_0)]\!]=[\![(\exists x\in z)\,\varphi(x)]\!]$.

{\bf(3)}
Let
$\Phi$
be a~correspondence from~$X$ to~$Y$
inside
${\mathbbm V}^{(B)}$.
In other words,
$\Phi$,
$X$, and~$Y$ are
members of~${\mathbbm V}^{(B)}$,
while
$[\![\Phi\subset X\times Y]\!]=\mathbb 1$.
Then there is a unique correspondence
~$\Phi{\downarrow}$
from~$X{\downarrow}$
to~$Y{\downarrow}$
such that
$
\Phi{\downarrow}(A{\downarrow})=\Phi(A){\downarrow}
$
for every nonempty subset~$A$
of~$X$
inside~${\mathbbm V}^{(B)}$.
The resulting correspondence~$\Phi{\downarrow}$
from
$X{\downarrow}$
to~$Y{\downarrow}$
is   the
{\it descent\/}
of~$\Phi$
from~$X$ to~$Y$
inside~${\mathbbm V}^{(B)}$.

{\bf(4)}
The descent of the composite of correspondences
inside~${\mathbbm V}^{(B)}$
is the composite of their descents:
$
(\Psi\circ\Phi){\downarrow}=\Psi{\downarrow}\circ\Phi{\downarrow}.
$

{\bf(5)}
If
$\Phi$
is a correspondence inside~${\mathbbm V}^{(B)}$
then
$
(\Phi^{-1}){\downarrow}=(\Phi{\downarrow})^{-1}.
$

{\bf(6)}
Let
${\rm Id}_X$
be the identity mapping inside~${\mathbbm V}^{(B)}$
of some~$X\in {\mathbbm V}^{(B)}$.
Then
$
({\rm Id}_X){\downarrow}={\rm Id}_{X{\downarrow}}.
$

{\bf(7)}
Assume that
$[\![f:X\to Y]\!]=\mathbb 1$,
with
$X,Y,f\in {\mathbbm V}^{(B)}$;
i.e.,
$f$
is a~funstion from~$X$
to~$Y$
inside~${\mathbbm V}^{(B)}$.
Then
$f{\downarrow}$
is the~unique mapping from~$X{\downarrow}$
to~$Y{\downarrow}$
such that
$[\![f{\downarrow}(x)=f(x)]\!]=\mathbb 1$ for all $x\in X{\downarrow}$.

By~\hbox{(1)--(7)}, we may view
the descent
as a functor from the category of
$B$-valued sets and mappings (correspondences)
to the usual category of  sets and mappings
(correspondences) in the sense of~${\mathbbm V}$.

{\bf(8)}
Given
$x_1,\dots,x_n\in {\mathbbm V}^{(B)}$,
denote by
$(x_1,\dots,x_n)^B$
the corresponding  $n$-tuple inside~${\mathbbm V}^{(B)}$.
Assume that
$P$ is an~$n$-ary relation on~$X$
inside~${\mathbbm V}^{(B)}$;
i.e.,
$X,P\in {\mathbbm V}^{(B)}$
and
$[\![P\subset X^n]\!]=\mathbb 1$,
with $n\in \omega$.
Then there is a~$n$-ary relation~$P'$
on~$X{\downarrow}$
such that
$$
(x_1,\dots,x_n)\in P'\leftrightarrow [\![(x_1,\dots,x_n)^B\in P]\!]=\mathbb 1.
$$
Slightly abusing notation, we will
denote~$P'$ by the already-occupied symbol~$P{\downarrow}$
and call $P{\downarrow}$ the {\it descent\/} of~$P$.
We descend functions of several variables by analogy.

\subsection{}{\bf Ascent.}
Assume that
$x\in {\mathbbm V}$
and
$x\subset {\mathbbm V}^{(B)}$;
i.e., let
$x$ be some set composed of
$B$-valued sets or, in other words,
$x\in \mathscr P({\mathbbm V}^{(B)})$.
Put
$\varnothing{\uparrow}:=\varnothing$,
while
$\dom (x{\uparrow}):=x$ and
$\im (x{\uparrow}):=\{\mathbb 1\}$
in case $x\neq \varnothing$.
The member
$x{\uparrow}$
of the (separated) universe~${\mathbbm V}^{(B)}$;
i.e., the distinguished representative
of the class
$\{y\in {\mathbbm V}^{(B)}\mid [\![y=~x{\uparrow}]\!]=\mathbb 1
 \}$,
is   the {\it ascent\/} of~$x$.

{\bf(1)}
For all
$x\in \mathscr P({\mathbbm V}^{(B)})$
and every formula~$\varphi$, the following are valid:
$$
%\gathered
[\![(\forall z\in x{\uparrow})\,\varphi(z)]\!]=
\bigwedge_{y\in x} [\![\varphi (y)]\!], \quad
%\\
[\![(\exists z\in x{\uparrow})\,\varphi(z)]\!]=
\bigvee_{y\in x} [\![\varphi (y)]\!].
%\endgathered
$$

Ascending the correspondence $\Phi\subset X\times Y$, we should keep in mind that
the domain of departure of $\Phi$, which is $X$,
 may differ from the domain of~$\Phi$, which is
 $\dom (\Phi):=\{x\in X\mid\Phi(x)\ne \varnothing \}$.
This circumstance is immaterial for our nearest goals. Therefore,
 speaking about ascents, we will always imply the
{\it total\/} correspondences $\Phi$ satisfying
$\dom (\Phi)=X$.

{\bf(2)}
Take
$X,Y,\Phi\in {\mathbbm V}$. Assume that
$X,Y\subset {\mathbbm V}^{(B)}$,
and
$\Phi$
is a~correspondence from~$X$
to~$Y$.
There is a unique correspondence~$\Phi{\uparrow}$
from~$X{\uparrow}$
to~$Y{\uparrow}$
inside~${\mathbbm V}^{(B)}$
such that  the equality
$
\Phi{\uparrow}(A{\uparrow})=\Phi(A){\uparrow}
$
holds for every subset~$A$
of~$\dom (\Phi)$
if and only if
$\Phi$ is
{\it extensional}; i.e.,
$\Phi$ enjoys the property
$$
y_1\in \Phi(x_1)\to [\![x_1=x_2]\!]\le
\bigvee_{y_2\in \Phi(x_2)} [\![y_1=y_2]\!]
$$
for
$x_1,x_2\in \dom (\Phi)$.
Moreover,
$\Phi{\uparrow}=\Phi'{\uparrow}$,
with
$\Phi':=\{(x,y)^B\mid (x,y)\in \Phi \}$.
The member
$\Phi{\uparrow}$
of~${\mathbbm V}^{(B)}$
is   the {\it ascent\/}
of~$\Phi$.

{\bf(3)}
The composite of extensional correspondences is extensional too.
Moreover, the ascent of the composite is the composite of the ascents
inside~${\mathbbm V}^{(B)}$;
i.e., if
$\dom (\Psi)\supset \im (\Phi)$
then
$
{\mathbbm V}^{(B)}\models(\Psi\circ\Phi){\uparrow}=\Psi{\uparrow}\circ\Phi{\uparrow}.
$

Note also that if
$\Phi$
and
$\Phi^{-1}$
are extensional then
$(\Phi{\uparrow})^{-1}=(\Phi^{-1}){\uparrow}$.
However, the extensionality of~$\Phi$
does not guarantee the extensionality of~$\Phi^{-1}$ in general.

{\bf(4)}
It is worth mentioning that
if some extensional correspondence~$f$ is a function from~$X$
to~$Y$ then the ascent~$f{\uparrow}$ of~$f$
will be a function from~$X{\uparrow}$
to~$Y{\uparrow}$.
In this event, we may write down the extensionality property as follows:
$[\![x_1=x_2]\!]\le [\![f(x_1)=f(x_2)]\!]$ for all $x_1,x_2\in X$.

Given
$X\subset {\mathbbm V}^{(B)}$,
we let the symbol~$\mix(X)$
stand for the collection of all mixings
like~$\mix (b_{\xi}x_{\xi})$,
where
$(x_{\xi})\subset X$
and
$(b_{\xi})$
is an~arbitrary partition of unity.
The propositions to follow are called
the  {\it arrow cancellation rules\/}
or the {\it rules for ascending and descending}.
There are many nice reasons to refer to them simply as
the  {\it Escher rules}~\cite{Hofstadter}.

{\bf(5)}
Let
$X$
and
$X'$
be two subsets of~${\mathbbm V}^{(B)}$,
and
let
$f:X\to X'$
be extensional. Assume that
$Y,Y',g\in {\mathbbm V}^{(B)}$
enjoy the property
$[\![\,Y\ne \varnothing]\!]=[\![\,g:Y\to Y']\!]=\mathbb 1$.
Then
$$
\gathered
X{\uparrow}{\downarrow}=\mix (X),\quad
Y{\downarrow}{\uparrow}=Y;\quad
\\
f{\uparrow}{\downarrow}|_{X}=f,\quad
g{\downarrow}{\uparrow}=g.
\endgathered
$$

{\bf(6)}~By analogy with 4.1\,(6), we easily infer the following convenient equality:
$$
\mathscr P_{\Fin}(X\!\uparrow)=\{{\theta\!\uparrow} \mid\theta\in\mathscr P_{\Fin}(X)\}\!\uparrow.
$$
We may generally assert inside~${\mathbbm V}^{(B)}$ only that
$\mathscr P(X\!\uparrow)\supset\{{\theta\!\uparrow} \mid
\theta\in\mathscr P(X)\}\!\uparrow$.

\subsection{}{\bf Modified descents and ascents.}
Let
$X\in {\mathbbm V}$ and
$X\ne\varnothing$;
i.e.,
$X$
is a~nonempty set. Denote by~$\imath$
the canonical embedding
$x\mapsto x^{\scriptscriptstyle\wedge}$
for $x\in X$.
Then
$\imath(X){\uparrow}=X^{\scriptscriptstyle\wedge}$
and
$X=\imath^{-1}(X^{\scriptscriptstyle\wedge}{\downarrow})$.
Using these equalities, we can translate the operations of
ascending and descending to the case in which
$\Phi$
is a~correspondence from~$X$ to~$Y{\downarrow}$
and
$[\![\Psi$ is a~correspondence from~$X^{\scriptscriptstyle\wedge}$
to~$Y]\!]=\mathbb 1$,
where
$Y\in {\mathbbm V}^{(B)}$.
Namely, we put
$\Phi\upwardarrow:=(\Phi\circ\imath^{-1}){\uparrow}$
and
$\Psi\downwardarrow:=\Psi{\downarrow}\circ\imath$.
In this event, the member
$\Phi\upwardarrow$
is   the
{\it modified ascent\/}
of~$\Phi$,
while
$\Psi\downwardarrow$
is   the
{\it modified descent\/}
of~$\Psi$.
If the context prevents us from ambiguity  then
we will continue to speak simply about ascents and descents, using
the usual arrows.
Clearly,
$\Phi\upwardarrow$
is the unique correspondence inside~${\mathbbm V}^{(B)}$
enjoying the equality
$
[\![\Phi\upwardarrow(x^{\scriptscriptstyle\wedge})=\Phi(x){\uparrow}]\!]=\mathbb 1
$
for all $x\in X$.
By analogy,
$\Psi\downwardarrow$
is the unique correspondence from~$X$
to~$Y{\downarrow}$ such that the equality
$
\Psi\downwardarrow(x)=\Psi(x^{\scriptscriptstyle\wedge}){\downarrow}
$
holds for all $x\in X$.
If
$\Phi:=f$
and~$\Psi:=g$
are functions then the above equalities take the shape
$
[\![f\upwardarrow(x^{\scriptscriptstyle\wedge})=f(x)]\!]=\mathbb 1$
and
$g\downwardarrow(x)=g(x^{\scriptscriptstyle\wedge})$
for all $x\in X$.

\section{\bf Elements of Boolean Valued Analysis}

Each Boolean valued universe has a complete mathematical toolkit
containing the sets with arbitrary additional structures:
groups, rings, algebras, etc.
Descending the algebraic systems from
Boolean valued models yields new objects with new properties,
revealing many facts about their constructions and interrelations.
Search into these interrelations is the content of Boolean valued
analysis. In this survey we naturally confine exposition only to the basics of
the relevant technique.

\subsection{}{\bf Boolean sets.}
We start with recognizing the simplest objects
that descend from a Boolean valued universe.

{\bf(1)}
A~{\it Boolean set},
or~a~{\it set with $B$-structure},
or simply a~{\it $B$-set} is by definition
a~pair $(X,d)$,
where
$X\in {\mathbbm V}$,
$X\ne \varnothing$,
and
$d$ is a~mapping from~$X\times X$
to a~Boolean algebra~$B$
satisfying for all $x,y,z\in X$ the following conditions:

\begin{itemize}
\item[\bf(a)]
$d (x,y)=\mathbb 0\leftrightarrow x=y$;

\item[\bf(b)]
$d (x,y)=d (y,x)$;

\item[\bf(c)]
$d (x,y)\le d (x,z)\vee d (z,y)$.
\end{itemize}

\noindent
Each nonempty subset
$\varnothing\ne X\subset {\mathbbm V}^{(B)}$
provides an example of a~$B$-set on assuming that
$d (x,y):=[\![x\ne y]\!]=\neg [\![x=y]\!]$ for all $x,\ y\in X$.
Another example arises if we furnish a
nonempty set~$X$
with the ``discrete $B$-metric''~$d$;
i.e., on letting
$d (x,y)=\mathbb 1$
in case
$x\ne y$
and
$d (x,y)=\mathbb 0$
in case
$x=y$.

{\bf(2)}
Given an arbitrary Boolean algebra $D$, we may take as a~$D$-metric
the {\it symmetric difference}:
$x\vartriangle y\!:= (x\wedge \neg\kern1pt y)\vee (y\wedge \neg\kern1pt x)$.
We  proceed  now with a relevant construction.
Let $\psi$ be an~ultrafilter on a~Boolean algebra~$D$.
Consider a~Boolean set $(X,d_X)$ with some $D$-metric~$d_X$.
Introduce the binary relation ${\sim_\psi}$ on~$X$ by the rule
 $
 (x,y)\in{\sim_\psi}\leftrightarrow \neg\, d_X(x,y)\in\psi.
 $
From the definition of Boolean metric it follows
that ${\sim_\psi}$ is an~equivalence.
Let $X/{\sim_\psi}$ stand for the factor-set of~$X$
by~${\sim_\psi}$, and let $\pi_X:X\to X/{\sim_\psi}$ be
the canonical factor-mapping.
If we start with the Boolean set $(D,\vartriangle)$
then we see that $D/{\sim_\psi}$ is the two-point
Boolean algebra~$\mathbb 2$.
In this event there is a unique mapping
$\overline{d}:X/{\sim_\psi}\to \mathbb 2$ such that
 $\overline{d}(\pi_X(x),\pi_X(y))=\pi_D(d(x,y))$ for all $x,y\in X$.
Moreover, $\overline{d}$ is the discrete  Boolean metric on~$X/{\sim_\psi}$.
If $d_X$ is the discrete metric on$X$ then
${\sim_\psi} ={\rm Id_X}$ and $X/{\sim_\psi}=X$.

{\bf(3)}
Let
$(X,d)$
be some $B$-set.
Assume that $\psi\!:=\psi_{{\rm Id}_B}$ is an ultrafilter
constructed on~$B^{\scriptscriptstyle\wedge}$ as in 4.1\,(7).
By the restricted transfer principle
 $(X^{\scriptscriptstyle\wedge},d^{\scriptscriptstyle\wedge})$
is a~$B^{\scriptscriptstyle\wedge}$-set inside~${\mathbbm V}^{(B)}$.
Put ${\mathscr X}:=X^{\scriptscriptstyle\wedge}/{\sim_\psi}$. In this event
 $[\![x^{\scriptscriptstyle\wedge} \sim_\psi y^{\scriptscriptstyle\wedge}]\!]
 =\neg\kern1pt d(x,y)$ for all $x,y\in X$.
Consequently, there are a~member
$\mathscr X$ of~${\mathbbm V}^{(B)}$
and an injection~$\imath:X\to X':=\mathscr X{\downarrow}$
such that
$d (x,y)=[\![\imath(x)\ne \imath(y)]\!]$
for all $x, y\in X$
and every~$x'\in X'$
admits the representation
$x'=\mix_{\xi\in \Xi}(b_{\xi}\imath (x_{\xi}))$,
with
$(x_{\xi})_{\xi\in \Xi}\subset X$
and
$(b_{\xi})_{\xi\in \Xi}$
a~partition of unity in~$B$.
The element~$\mathscr X$ of~${\mathbbm V}^{(B)}$
is   the {\it Boolean valued realization\/}
of the $B$-set~$X$.
If
$X$
is a discrete $B$-set then
$\mathscr X=X^{\scriptscriptstyle\wedge}$
and
$\imath(x)=x^{\scriptscriptstyle\wedge}$
for all $x\in X$.
If
$X\subset {\mathbbm V}^{(B)}$
then
$\imath{\uparrow}$
is an injection from~$X{\uparrow}$
to~$\mathscr X$
inside~${\mathbbm V}^{(B)}$.

{\bf(4)}
A mapping~$f$
from a~$B$-set
$(X,d)$
to a~$B$-set
$(X',d\kern1pt')$
is   {\it contractive\/}
provided that
$d\kern1pt'(f(x),f(y))\le d(x,y)$
for all
$x,y\in X$.

Assume that
$X$
and
$Y$
are some $B$-sets. Assume further that
$\mathscr X$
and
$\mathscr Y$
are the Boolean valued realizations of~$X$ and~$Y$,
while
$\imath:X\to \mathscr X{\downarrow}$
and
$\jmath:Y\to \mathscr Y{\downarrow}$
are the corresponding injections.
If
$f:X\to Y$
is a contractive mapping then there is a unique
member~$g$ of~${\mathbbm V}^{(B)}$
such that
$[\![g:\mathscr X\to \mathscr Y]\!]=\mathbb 1$
and
$f=\jmath^{-1}\circ g{\downarrow}\circ\imath$.
We agree to let
$\mathscr X:=\mathscr F^\sim(X):=X^\sim$
and
$g:=\mathscr F^\sim(f):=f^\sim$.

{\bf(5)}
The following hold:

\begin{itemize}
 \item[\bf(a)]
 {${\mathbbm V} ^{(B)}\models f{(A)^\sim}={f^\sim}({A^\sim} )$
for $A\subset X$};

\item[\bf(b)]
If $g:Y\to Z$ is a contractive mapping then
so is  $g\circ f$ and, moreover,
${\mathbbm V} ^{(B)}\models(g\circ f)^\sim ={g^\sim}\circ{f^\sim} $;

\item[\bf(c)]
{${\mathbbm V} ^{(B)}\models
``{f^\sim}$ is an injection'' if and only if
 $f$ is a~$B$-isometry};

\item[\bf(d)]
{${\mathbbm V} ^{(B)}\models $
``${f^\sim}$ is a surjection'' if and only if
 $$(\forall y\in Y)\ \bigvee \{d(f(x),y)\mid x\in X\}=\mathbb 1.$$}
\end{itemize}

{\bf(6)}
 Consider a $B$-set $(X,d)$.
Let $(b_\xi)$ be a~partition of unity of~$B$
and let $(x_\xi)$ be a~family of elements of~$X$.
A~{\it mixing\/} of~$(x_\xi)$ by~$(b_\xi)$
 is  $x\in X$ such that $b_\xi\wedge d(x,x_\xi)={\mathbb 0}$ for all~$\xi$.
We denote the instance of mixing as before: $x=\mix (b_\xi x_\xi)$.
 If a mixing exists then it is unique.
 Indeed, if
 $y\in X$ and $(\forall\,\xi) (b_\xi\wedge d(y,x_\xi)={\mathbb 0})$
 then
 $$
 b_\xi\wedge d(x,y)\le b_\xi\wedge (d(x,x_\xi)\vee d(x_\xi,y))={\mathbb 0}.
 $$
 Every complete Boolean algebra $B$ enjoys the infinite distributive law
 and so
 $$
 d(x,y) =\bigvee\{b_\xi\wedge d(x,y)\}={\mathbb 0}.
 $$
 Hence, $x=y$.

 Note that, in contrast to the case of~${\mathbbm V}^{(B)}$ (cp.~4.3),
 some mixings may fail to exist in an~abstract $B$-set.

{\bf(7)}
 Consider a~$B$-set $(X,d)$.
 Given $A\subset X$, denote by $\mix (A)$ the set of all mixings
 of the  members of~$A$. If $\mix (A)=A$ then we say that $A$
 is a~{\it cyclic subset\/} of~$X$.
 The intersection of all cyclic sets  including $A$
 is denoted by $\cyc(A)$.
 A Boolean set $X$ is  {\it universally complete\/}
or {\it extended\/}
provided that $X$ contains all mixings  $\mix (b_\xi x_\xi)$
of arbitrary families $(x_\xi)\subset X$ by all partitions of unity $(b_\xi)\subset B$.
In case these mixings exist for
finite families,  $X$ is  {\it finitely complete\/} or {\it decomposable}.
It can be shown that if $X$ is a~universally complete
 $B$-set then $\mix (A)=\cyc(A)$ for every $A\subset X$.
A cyclic subset of a~$B$-set may fail to be universally complete.
At the same time, each cyclic subset of~${\mathbbm V}^{(B)}$, furnished with
the canonical $B$-metric, is a~universally complete $B$-set.

\subsection{}{\bf Algebraic  ${\boldsymbol B}$-systems.}
Recall that by a~{\it signature\/}
we mean the 3-tuple
$\sigma:=(F,P,{\mathfrak a})$, where
$F$ and $P$  are some (possibly empty) sets,
while
${\mathfrak a}$ is a~mapping from~$F\cup P$ to~$\omega$.
If $F$ and $P$ are finite then  $\sigma$
is   a~{\it finite signature}.
In applications we usually deal with
the algebraic systems of finite signature.
We thus confine exposition to considering only finite signatures.

By an~$n$-{\it ary operation\/}
and
 {\it $n$-ary predicate\/}
on a~$B$-set
$A$ we mean some contractive mappings $f:A^n\to A$
and $p:A^n\to B$.
Recall that  $f$ and $p$
are {\it  contractive\/}
provided that
$$
\gathered
d(f(a_0,\dots,a_{n-1}),f(a'_0,\dots,a'_{n-1}))\le
\bigvee\limits_{k=0}^{n-1}\,d(a_k,a'_k),
\\
\vartriangle\!\big(p(a_0,\dots,a_{n-1}),
p(a'_0,\dots,a'_{n-1})\big)
\le\bigvee\limits_{k=0}^{n-1}\,
d(a_k,a'_k)
\endgathered
$$
for all $a_0$, $a'_0,\dots,a_{n-1}$,
$a'_{n-1}\in A$, where $d$ is the~$B$-metric on~$A$.

The above definitions depend explicitly on~$B$ and so it would
be wise to speak of $B$-operations, $B$-predicates, etc.
However, it is wiser to simplify and economize
whenever this involves no confusion.

An {\it algebraic $B$-system\/}
${\mathfrak A}$
of signature~$\sigma$
is by definition a~pair $(A,\nu)$,
where $A$ is a~nonempty
$B$-set,
the {\it carrier},
or
{\it underlying set},
or
{\it universe\/}
of~${\mathfrak A}$,
while $\nu$ is a~mapping such that
\begin{itemize}
\item[\bf(a)]
$\dom(\nu)=F\cup P$;
\item[\bf(b)]
$\nu (f)$ is an~${\mathfrak a}(f)$-ary operation on~$A$
for all $f\in F$;
\item[\bf(c)]
$\nu (p)$  is an~${\mathfrak a}(p)$-ary predicate on~$A$
for all $p\in P$.
\end{itemize}
\noindent
We call $\nu$
the {\it interpretation\/} of~${\mathfrak A}$
and write  $f^\nu$ and $p^{\kern1pt\nu}$
instead of $\nu (f)$ and $\nu (p)$.
The signature of an~algebraic $B$-system ${\mathfrak A}:=(A,\nu)$
is often denoted by~$\sigma({\mathfrak A})$
and the carrier $A$ of~${\mathfrak A}$, by~$|{\mathfrak A}|$.

An algebraic $B$-system $\mathfrak A:=(A,\nu)$
is {\it universally complete\/} or
{\it finitely complete\/} provided that
so is the carrier of~$A$.
Since  $A^0 =\{\varnothing\}$, the  0-ary operations
and predicates on~$A$ are mappings from~$\{\varnothing\}$
to~$A$ and~$B$ respectively.
We will identify the mapping $g:\{\varnothing\}\to A\cup B$
with the set $g(\varnothing)$.
Each 0-ary operation on~$A$ transforms in this fashion into
a~unigue member of~$A$. Similarly, the set of all
0-ary predicates on~$A$ becomes a subset of the Boolean
algebra~$B$.

If $F:=\{f_1,\dots,f_n\}$ and
$P:=\{p_1,\dots,p_m\}$ then an algebraic $B$-system
of signature $\sigma$ is written  as
$(A,\nu (f_1),\dots,\nu (f_n)$,
$\nu (p_1),\dots,\nu (p_m))$,
or even
$(A,f_1,\dots,f_n$, $p_1,\dots,p_m)$.
Moreover, we will substitute
$\sigma=(f_1,\dots,f_n$,
$p_1,\dots,p_m)$
 for~$\sigma=(F,P,{\mathfrak a})$.

We  turn now to the $B$-valued interpretation of
a first order language. Consider an algebraic
$B$-system ${\mathfrak A}:=(A,\nu)$ of signature
$\sigma:=\sigma ({\mathfrak A}):=(F,P,{\mathfrak a})$.

Let $\varphi(x_0,\dots,x_{n-1})$ be some formula of
signature~$\sigma$ with~$n$  free variables. Assume given
$a_0,\dots,a_{n-1}\in A$.
In this event we may define the Boolean truth value
$$
|\varphi |^{{\mathfrak A}}
(a_0,\dots,a_{n-1})\in B
$$
of a formula $\varphi$ in~${\mathfrak A}$ for the given values
 $a_0,\dots,a_{n-1}$ of the variables $x_0,\dots,x_{n-1}$.
The definition proceeds by usual recursion on the length
of~$\varphi$.
Considering the logical connectives and quantifiers,
put
$$
\gathered
|\varphi\wedge\psi |^{{\mathfrak A}}\, (a_0,\dots,a_{n-1}):=
|\varphi |^{\mathfrak A} (a_0,\dots,a_{n-1})\wedge |\psi |^{{\mathfrak A}}
(a_0,\dots,a_{n-1});
\\
|\varphi\vee\psi |^{\mathfrak A}\,(a_0,\dots,a_{n-1}):=
|\varphi |^{\mathfrak A} (a_0,\dots,a_{n-1})\vee |\psi |^{{\mathfrak A}}
(a_0,\dots,a_{n-1});
\\
|\neg\varphi |^{{\mathfrak A}}\,(a_0,\dots,a_{n-1}):=\neg |\varphi |^{{\mathfrak A}} (a_0,
\dots,a_{n-1});
\\
|(\forall\,x_0)\varphi |^{{\mathfrak A}}\, (a_1,\dots,a_{n-1}):=
\bigwedge\limits_{a_0\in A}|\varphi |^{{\mathfrak A}} (a_0,\dots,a_{n-1});
\\
|(\exists\,x_0)\varphi |^{{\mathfrak A}}\,(a_1,\dots,a_{n-1}):=
\bigvee\limits_{a_0\in A}|\varphi |^{{\mathfrak A}} (a_0,\dots,a_{n-1}).
\endgathered
$$
We turn now to the atomic formulas.
Suppose that  $p\in P$ stands for an~$m$-ary predicate,
$q\in P$  is a~0-ary predicate, and
$t_0,\dots,t_{m-1}$ are some terms of signature~$\sigma$
that take the values $b_0,\dots,b_{m-1}$
when the variables $x_0,\dots,x_{n-1}$
assume the values $a_0,\dots,a_{n-1}$.
We proceed with letting
$$
\gathered
|\varphi |^{{\mathfrak A}} (a_0,\dots,a_{n-1}):=\nu (q),
\mbox{~if~}\varphi=q;
\\
|\varphi |^{{\mathfrak A}} (a_0,\dots,a_{n-1}):=\neg d(b_0,b_1),
\,\kern5pt\mbox{~if~}\varphi=(t_0 =t_1);
\\
|\varphi |^{{\mathfrak A}} (a_0,\dots,a_{n-1}):=p^\nu
(b_0,\dots,b_{m-1}),
\,\kern5pt\mbox{~if~}\varphi=p (t_0,\dots,t_{m-1}),
\endgathered
$$
where $d$ is  the~$B$-metric on~$A$.

A formula
$\varphi(x_0,\dots,x_{n-1})$
is {\it satisfied\/}
in an algebraic $B$-system
${{\mathfrak A}}$ by the assignment
$a_0,\dots,a_{n-1}\in A$
of  $x_0,\dots,x_{n-1}$,
in symbols ${{\mathfrak A}}\models\varphi(a_0,\dots,a_{n-1})$,
provided that
$|\varphi |^{{\mathfrak A}} (a_0,\dots,a_{n-1})=\mathbb 1_B$.
The concurrent expressions are as follows:
$a_0,\dots,a_{n-1}\in A$
{\it satisfy\/}
$\varphi(x_0,\dots,x_{n-1})$
or $\varphi(a_0,\dots,a_{n-1})$
{\it holds\/} in~${{\mathfrak A}}$.
In case of the two-point Boolean algebra, we come to the conventional
definition of satisfaction in an~algebraic system.

Recall that a closed formula
$\varphi $
of signature $\sigma$  is a~{\it tautology\/}
or {\it logically valid formula\/}
provided that
$\varphi$ holds in every algebraic
$\mathbb 2$-system of signature~$\sigma$.

{\bf(1)}
Consider some algebraic $B$-systems
${{\mathfrak A}}:=(A,\nu)$
and $\mathfrak D:=(D,\mu)$ of the same signature~$\sigma$.
A~mapping $h:A\to D$ is a~{\it homomorphism\/}
from~${{\mathfrak A}}$ to~$\mathfrak D$
provided that the following hold for all
$a_0,\dots,a_{n-1}\in A$:

\begin{itemize}
\item[\bf(a)]
$d_D(h(a_1),h(a_2))\le d_A(a_1,a_2)$;

\item[\bf(b)]
$h(f^\nu)=f^\mu$ if
${{\mathfrak a}}(f)=0$;

\item[\bf(c)]
$h(f^\nu (a_0,\dots,a_{n-1}))=
f^\mu (h(a_0),\dots,h(a_{n-1}))$ if
$0\ne n :={{\mathfrak a}}(f)$;

\item[\bf(d)]
$p^{\,\nu} (a_0,\dots,a_{n-1})\le p^{\,\mu} (h(a_0),\dots,h(a_{n-1}))$
where $n :={{\mathfrak a}}(p)$.

\end{itemize}
\medskip

{\bf(2)}
A homomorphism
$h$ is {\it strong\/}
provided that
for  $p\in P$
such that
$0\ne n:=\mathfrak a (p)$
we have
$$
\gathered
p^{\,\mu} (d_0,\dots,d_{n-1})
\\
\le\,\bigvee\limits_{a_0,\dots,a_{n-1}
\in A}\,{p^{\,\nu} (a_0,\dots,a_{n-1})\wedge
\neg d_D(d_0,h(a_0))\wedge
\dots\wedge
\neg d_D(d_{n-1},h(a_{n-1}))}
\endgathered
$$
\noindent
for all $d_0,\dots,d_{n-1}\in D$.

If  (a) and (d)
turn into equalities then
$h$ is an~{\it isomorphism\/} from~$\mathfrak A$ to~$\mathfrak D$.
Obviously, each surjective isomorphism (in particular,
the identity mapping
${\rm Id}_A:A\to A$)
is a strong homomorphism.
The composite of (strong) homomorphisms is
a~(strong) homomorphism as well.
Clearly, if $h$ is a~homomorphism and
$h^{-1}$  is a~homomorphism too,
then  $h$ is an~isomorphism.
We mention again that in the case of the two-point
Boolean algebra we come to the conventional concepts
of homomorphism, strong homomorphism, and isomorphism.

\subsection{}{\bf  Descending~2.}
Before defining the descent of a~general
algebraic system,  consider the descent of the two-point
Boolean algebra.
Take two arbitrary  members
$0$, $1\in{\mathbbm V} ^{(B)}$
satisfying the condition $[\![0\ne 1]\!]=\mathbb 1_B$.
For instance, put $0:=\mathbb 0_B^{\scriptscriptstyle\wedge}$
and
$1:=\mathbb 1^{\scriptscriptstyle\wedge} _B$.
The descent $D$ of the two-point Boolean algebra~$\{0,1\}^B$
inside~${\mathbbm V} ^{(B)}$
is a~complete Boolean algebra isomorphic to~$B$.
The formulas
$
[\![\chi(b)=1]\!]=b$ and
$[\![\chi (b)=0]\!]=\neg b$, with $b\in B$,
yield the isomorphism $\chi:B\to D$.

\subsection{}{\bf  Descending an algebraic $\boldsymbol B$-system.}
Look now at  an~algebraic system ${{\mathfrak A}}$ of
signature~$\sigma$ inside~${\mathbbm V} ^{(B)}$
and assume that
${\mathbbm V} ^{(B)}\models$ ``${{\mathfrak A}} =(A,\nu)$
for some $A$ and~$\nu$.''
The {\it descent\/} of~${{\mathfrak A}}$ is the pair
${{{\mathfrak A}}\!\!\downarrow}:=(A\!\!\downarrow,\mu)$, where $\mu $
is the function defined as follows:
$f^\mu:=f^\nu{\downarrow}$ for  $f\in F$ and $p\kern1pt^\mu:=\chi^{-1}\circ p^{\kern1pt\nu}{\downarrow}$
for $p\in P$.
Here $\chi$ is
the isomorphism of~5.3 between the Boolean algebras
$B$ and~$\{0,1\}^B\!\!\downarrow$.

{\bf(1)}
Let $\varphi(x_0,\dots,x_{n-1})$ be a~distinguished formula
of signature~$\sigma$ with $n$ free variables.
Write down the formula $\Phi (x_0,\dots,x_{n-1},\!{{\mathfrak A}})$
of the language of set theory which expresses
the fact ${{\mathfrak A}}\models
\varphi(x_0$,$\dots,x_{n-1})$.
The formula
${{\mathfrak A}}\models\varphi(x_0,\dots,x_{n-1})$
determines an~$n$-ary predicate on~$A$ or, in other words,
a~mapping from $A^n$ to~$\{0,1\}$.
By the  maximum and transfer principles there is a~unique
member
$|\varphi |^{{\mathfrak A}}$ of~${\mathbbm V} ^{(B)}$
such that
$$
\gathered\,
[\![|\varphi |^{{\mathfrak A}}:A^n\to\{0,1\}^B]\!]=\mathbb 1,
%\quad
\\
[\![|\varphi |^{{\mathfrak A}} (a_0,\dots,a_{n-1})=1]\!]=
[\![\Phi (a_0,\dots,a_{n-1},{{\mathfrak A}})]\!]=\mathbb 1
\endgathered
$$
for all  $a_0,\dots,a_{n-1}\in {A\!\!\downarrow}$.
Therefore, the formula
$$
{\mathbbm V} ^{(B)}\models \mbox{\rm ``}\varphi(a_0,\dots,a_{n-1})
\mbox{\rm\ holds\  in~}{\mathfrak A}\mbox{\rm ''}
$$
is valid if and only if
$[\![\Phi (a_0,\dots,a_{n-1},{{\mathfrak A}})]\!]=\mathbb 1$.

{\bf(2)}
Let ${{\mathfrak A}}$ be an~algebraic system
of signature~$\sigma$
inside~${\mathbbm V} ^{(B)}$.
Then ${{\mathfrak A}}\!\downarrow$ is a~universally complete
algebraic $B$-system of signature~$\sigma$.
Moreover,
$$
\chi\circ |\varphi |^{{{\mathfrak A}}\downarrow}={|\varphi |^{{\mathfrak A}}\!\!\downarrow}
$$
for every formula $\varphi$  of signature~$\sigma$.

{\bf(3)}
Let ${{\mathfrak A}}$ and $\mathfrak B$  be two algebraic systems
of the same signature~$\sigma$
inside~${\mathbbm V} ^{(B)}$.
Put ${{\mathfrak A}} ':={{\mathfrak A}}\!\!\downarrow$ and
$\mathfrak B ':=\mathfrak B\!\!\downarrow$.
Then if $h$ is a~homomorphism (strong homomorphism)
inside~${\mathbbm V} ^{(B)}$ from~${{\mathfrak A}}$ to~$\mathfrak B$
then
$h':=h\!\!\downarrow$ is a~homomorphism (respectively, strong homomorphism)
between the $B$-systems ${{\mathfrak A}} '$ and $\mathfrak B '$.

Conversely, if $h':{{\mathfrak A}} '\to\mathfrak B '$ is a~homomorphism
(strong homomorphism) of the relevant algebraic
$B$-systems then $h:=h'\!\!\uparrow $
is a~homomorphism  (respectively, strong homomorphism)
from~${{\mathfrak A}}$ to~$\mathfrak B$
inside~${\mathbbm V} ^{(B)}$.

{\bf(4)}
Let  ${\mathfrak A}:=(A,\nu)$  be an~algebraic $B$-system of signature~$\sigma$.
Then there are  $\mathscr A$ and $\mu\in{\mathbbm V} ^{(B)}$ such that

\begin{itemize}
\item[\bf(a)]
${\mathbbm V} ^{(B)}\models $
``$(\mathscr A,\mu)$ is an~algebraic system of
signature~$\sigma $''$;$

\item[\bf(b)]
if ${{\mathfrak A}} ':=(A',\nu ')$  is the descent of~$(\mathscr A,\mu)$
then ${{\mathfrak A}}$ is a~universally complete
algebraic $B$-system of signature~$\sigma;$

\item[\bf(c)]
there is an isomorphism $\imath$ from~${{\mathfrak A}}$ to~${{\mathfrak A}}'$
such that $A'=\mix(\imath(A));$

\item[\bf(d)]

given a formula~$\varphi$ of signature $\sigma$ in  $n$ free variables,
the following holds
$$
\gathered
|\varphi |^{{\mathfrak A}} (a_0,\dots,a_{n-1})=|\varphi |^{{{\mathfrak A}} '}
(\imath(a_0),\dots,\imath(a_{n-1}))
\\
=\chi ^{-1}\bigl(
(|\varphi |^{{{\mathfrak A}}^\sim})\!\!\downarrow\!\!(\imath(a_0),\dots,
\imath(a_{n-1}))\bigr)
\endgathered
$$
for all $a_0,\dots,a_{n-1}\in A$, where $\chi$ is defined as above.
\end{itemize}

\section{\bf Boolean Valued Reals}
We are now able to apply the technique of
Boolean valued analysis to the algebraic system
of the utmost importance for the whole body of mathematics,
the field of the reals.
Let us recall a few preliminaries
of the theory of vector lattices.

\subsection{}{\bf Kantorovich spaces.}
We start with some definitions.
Members $x$ and $y$ of a vector lattice $E$ are
{\it disjoint\/} (in symbols, $x\perp y$) provided that $|x| \wedge  |y| =0$.
A~{\it band\/} of~$E$ is the {\it disjoint complement\/}
$M^\perp\!:= \{x\in E \mid (\forall y\in M)\, x\perp y\}$
of some  subset $M\subset E$. If this is sensible then
for simplicity we will exclude from consideration
the trivial vector lattice consisting of the sole zero.

The inclusion-ordered set ${\mathfrak B}(E)$ of the bands of~$E$
is a~complete Boolean algebra under the following operations:
$$
L\wedge K=L\cap K,\quad
L\vee K=(L\cup K)^{\perp\perp},\quad
\neg L =L^\perp\quad
(L,K\in{\mathfrak B}(E)).
$$
The Boolean algebra ${\mathfrak B}(E)$ is often referred to as
the {\it base\/} of~$E$.

A~{\it band projection\/} in~$E$ is an~idempotent linear operator
$\pi:E\to E$ such that $0\leq\pi x\leq x$
for all $0\leq x\in E$. The set ${\mathfrak P}(E)$ of all band projections,
ordered by the rule $\pi\le\rho\Longleftrightarrow\pi\circ\rho=\pi$,
is a~Boolean algebra with the following operations:
$$
\pi\wedge \rho =\pi\circ\rho,\quad
\pi\vee\rho =\pi +\rho -\pi\circ\rho,\quad
\neg\pi =I_E-\pi\quad (\pi,\rho\in\mathfrak (E)).
$$

Assume that $u\in E_+$ and $e\wedge (u-e)= 0$ for some $0\leq e\in E$.
Then $e$ is a~{\it fragment\/} or {\it component\/} of~$u$.
The set ${\mathfrak E}(u)$ of all fragments of~a~nonzero element $u$ with the induced
order from~$E$ is a~Boolean algebra with the same lattice operations as
in~$E$ and the complement acting by the rule $\neg e \!:= u - e$.

A vector lattice~$E$ is a~{\it Kantorovich space\/} or
$K$-{\it space} provided that $E$ is {\it Dedekind complete}; i.e.,
 each upper bounded nonempty subset of~$E$ has the least upper bound.

If each family of pairwise disjoint elements of a~$K$-space~$E$
is bounded then $E$ is  {\it universally complete\/} or
{\it extended}.

{\bf(1)}
Let  $E$ be an arbitrary \hbox{$K$-}space. Then the mapping
 $\pi\mapsto\pi(E)$ is an isomorphism between the Boolean
 algebras~${\mathfrak P}(E)$ and~${\mathfrak B}(E)$. If $E$
has an order unit ${\mathbb 1}$ then the mappings $\pi\mapsto\pi{\mathbb 1}$
from~${\mathfrak P}(E)$ to~${\mathfrak E}(\mathbb 1)$ and
$e\mapsto \{e\}^{\perp\perp}$ from~${\mathfrak E}(\mathbb 1)$ to~${\mathfrak B}(E)$
are isomorphisms between the corresponding
Boolean algebras.

{\bf(2)}
Each universally complete $K$-space $E$ with order unit~${\mathbb 1}$
admits a unique multiplication that
makes~$E$ into a faithful  $f$-algebra and ${\mathbb 1}$ into the ring unit.
Each band projection $\pi\in{\mathfrak P}(E)$
coincides with multiplication by~$\pi({\mathbb 1})$
in the resulting $f$-algebra.

\subsection{}{\bf  Descending the reals.}
By the maximum and transfer principles there is a~member $\mathscr R$ of~${\mathbbm V} ^{(B)}$
such that
${\mathbbm V} ^{(B)}\models $
``$\mathscr R$ is an ordered field of the reals.''
Clearly,
the field $\mathscr R$ is unique up to isomorphism inside ${\mathbbm V} ^{(B)}$;
i.e., if  ${\mathscr R}\kern1pt'$ is another field of
the reals inside~${\mathbbm V} ^{(B)}$ then
${\mathbbm V} ^{(B)}\models$
``the fields $\mathscr R$ and ${\mathscr R}\kern1pt'$ are isomorphic.''
It is easy that $\mathbbm R ^{\scriptscriptstyle\wedge} $
is an Archimedean ordered field inside~${\mathbbm V} ^{(B)}$.
Therefore we may assume that ${\mathbbm V} ^{(B)}\models $
``$\mathbbm R ^{\scriptscriptstyle\wedge}\subset\mathscr R$ and
$\mathscr R$ is a~(metric) completion of~$\mathbbm R ^{\scriptscriptstyle\wedge}$.''
Given the usual unity~$1$ of~$\mathbbm R$, note that
${\mathbbm V} ^{(B)}\models $
``$1:=1^{\scriptscriptstyle\wedge}$
is an order unit of~$\mathscr R$.''

Consider the descent $\mathscr R\!\!\downarrow$ of the algebraic system
$\mathscr R:=(|\mathscr R|,+,\,\cdot\,,0,1,\le)$.

{\bf(1)} In accord with the general construction,
we make the descent of the carrier of~$\mathscr R$
into an~algebraic system by descending the operations and order
on~$\mathscr R$.
In more detail, if we are given
$x,y,z\in\mathscr R\!\!\downarrow$ and
$\lambda\in\mathbbm R~$ then
we define the addition, multiplication, and order on $\mathscr R$
 as follows:

\begin{itemize}
\item[\bf(a)]$x+y=z\leftrightarrow [\![x+y=z]\!]=\mathbb 1;$

\item[\bf(b)]$xy=z\leftrightarrow [\![xy=z]\!]=\mathbb 1;$

\item[\bf(c)]$x\le y\leftrightarrow [\![x\le y]\!]=\mathbb 1;$

\item[\bf(d)]$\lambda x=y\leftrightarrow [\![\lambda ^{\scriptscriptstyle\wedge} x=y]\!]=\mathbb 1.$

\end{itemize}

{\bf(2)}
The algebraic system $\mathscr R\!\!\downarrow$ is a~universally complete
$K$-space. Moreover, there is a~(canonical)
isomorphism ~$\chi$
of the Boolean algebra~$B$
to the base ${\mathfrak P} ({\mathscr R\!\downarrow})$
of~$\mathscr R\!\!\downarrow$
such that
$$
\gathered
\chi (b)x=\chi (b)y\leftrightarrow b\le [\![x=y]\!],
\\[2\jot]
%\quad
\chi (b)x\le\chi (b)y\leftrightarrow b\le [\![x\le y]\!]
\endgathered
$$
for all $x,y\in\mathscr R\!\!\downarrow$ and $b\in B$.

This remarkable result, establishing the immanent connection
between  Boolean valued analysis and vector-lattice theory,
belongs to Gordon~\cite{Gordon}. The Gordon theorem has demonstrated
that each universally complete Kantorovich space provides
a~new model of the field of the reals, and all these models
have the same rights in mathematics.
Moreover, each Archimedean vector lattice, in particular,
an arbitrary  $L_p$-space, $p\ge 1$,  ascends into a dense
sublattice of the reals $\mathbbm R$ inside an appropriate
Boolean valued universe.

Descending the basic scalar fields opens  a turnpike
to the intensive application of  Boolean valued models in
functional analysis.
The technique of Boolean valued analysis demonstrates its efficiency
in studying Banach spaces and algebras as well as lattice-normed spaces and
modules. The corresponding results are collected and elaborated
in~\cite[Chs.~10--12]{IBA}.

\section{\bf The Cardinal Shift}

Boolean valued models were invented for
research into the foundations of mathematics.
Many delicate properties of the objects inside  ${\mathbbm V}^{(B)}$
depend essentially on the structure of the initial Boolean algebra~$B$.
The diversity of opportunities together with
a~great stock of information on
particular Boolean algebras ranks  Boolean valued models among the
most powerful tools of foundational studies.

The working mathematician feels rarely if ever
the dependence of his or her research on
foundations. Therefore, we allotted the bulk of the above exposition
to the constructions and properties of~$\mathbbm V^{(B)}$
in which the particularities of~$B$ were completely immaterial.
The basics of Boolean valued analysis rest precisely on these
constructions and properties. Neglect notwithstanding,
foundations underlie  the core of mathematics. Boolean valued analysis
stems from the brilliant results of
 G\"odel and Cohen who demonstrated to all of us the
 independence of the continuum hypothesis from the axioms
 of~$\ZFC$. It thus stands to reason to discuss the
 relevant matters in some detail.

\subsection{}{\bf Mysteries of the continuum.}
The concept of continuum   belongs to
the  most important general tools of science.
The mathematical views of the continuum
relate to the  understanding of time and time-dependent processes in physics.
It suffices to mention the great Newton and Leibniz who had
different perceptions of the continuum.

The smooth and perpetual motion, as well as vision of the nascent and
evanescent arguments
producing continuous changes in the dependant variables, underlies
Newton's outlook, philosophy, and his method of prime and ultimate ratios.
The principal difficulty of the views of Newton rests in the
impossibility of imagining the immediately preceding moment
of time nor the nearest neighbor of a given point of the continuum.
As regards Leibniz, he viewed every varying quantity as piecewise constant
to within higher order imperceptible infinitesimals.
His continuum splits into a collection of disjoint monads,
these immortal and mysterious ideal entities.

The views of Newton and Leibniz summarized the ideas that stem
from the remote ages. The mathematicians of Ancient Greece
distinguished between points and monads, so explicating the dual
nature of the objects of mathematics. The mystery of the structure
of the continuum came to us from our ancestors through two millennia.

The set-theoretic stance   revealed
a~new enigma of the continuum.
Cantor demonstrated that the set of the naturals is not equipollent
with the simplest mathematical continuum, the real axis.
This gave an immediate rise to the problem of the continuum
which consists in determining the cardinalities of the intermediate sets between
the naturals and the reals. The continuum hypothesis
reads that the intermediate subsets possess no new cardinalities.

The continuum problem was the first in the already cited
report by Hilbert~\cite{Problems}.
An incontrovertible {\it anti-ignorabimus}, Hilbert
was always inclined to the validity of the continuum hypothesis.
It is curious that one of his most beautiful and appealing
articles~\cite{Hilbert}, which is dated as of 1925  and contains
the famous phrase about Cantor's paradise, was devoted in fact
to an erroneous proof of the continuum hypothesis.

The Russian prophet  Luzin viewed as implausible even
the mere suggestion of the independence of the continuum hypothesis.
He said in his
famous talk ``Real Function Theory: State of the Art''~\cite{Luzin}
at the All-Russia Congress of Mathematicians
in~1927:
``The first idea that might leap to mind is that
the determination of the cardinality of the continuum is
a matter of a free axiom like the parallel postulate
of geometry. However, when we vary the parallel postulate,
keeping intact the rest of the axioms of Euclidean geometry,
we in fact change the precise meanings of the words
we write or utter, that is, `point,' `straight line,' etc.
What words are to change their meanings
if we attempt at making
the cardinality of the continuum movable along
the scale of alephs, while constantly proving consistency of this
movement?
The cardinality of the continuum, if only we imagine the latter as
a~set of points, is some unique entity that must reside in the
scale of alephs at the place which the cardinality of the continuum
belongs to;
no matter whether the determination of
this place is  difficult or even `impossible for us, the human beings' as
J.~Hadamard might comment.''

G\"odel proved the consistency of the continuum hypothesis
with the axioms of~$\ZFC$, by inventing the universe of
constructible sets~\cite{Goedel}.
Cohen demonstrated the consistency of the negation of the continuum hypothesis
with the axioms of $\ZFC$ by {\it forcing}, the new method
he invented for changing the properties of available or hypothetical
models of set theory.
 Boolean valued models made Cohen's difficult result
 simple,\footnote{Cohen remarked once about his result
 as follows \cite[p.~82]{Yandell}: ``But of course it {\it is\/} easy in the sense that there
 is a clear philosophical idea.''}
 demonstrating to the working mathematician the independence of the
 continuum hypothesis with the same visuality as the Poincar\'e model
 for non-Euclidean geometry.
Those who get acquaintance with this technique are inclined to follow
Cohen~\cite{Cohen} and view the continuum hypothesis
as ``{\it obviously} false.''

 We will keep to the famous direction of Dedekind
 that ``whatever is provable must not be believed without proof.''
 So we will briefly discuss  the formal side of the matter.
 More sophisticated arguments are collected in the remarkable
 articles by Kanamori~\cite{Kanamori},
 Cohen~\cite{Cohen-2002},  and  Manin~\cite{Manin}.

 \subsection{}{\bf Ordinals.}
 The concept of ordinal reflects the ancient trick of
 counting with successive notches.

A class $X$ is {\it transitive\/} provided that
each member of~$X$ is a~subset of~$X$; i.e.,
$$
 \Tr(X)\!:=(\forall\,y)
 (y\in X\to y\subset X).
 $$

An {\it ordinal class\/}
is a transitive class well ordered by the membership relation~$\in$.
The record $\Ord(X)$ means that $X$ is an~ordinal class.
An ordinal class presenting a~set is
an {\it ordinal\/} or a~{\it transfinite\/} number.
The class of all ordinals is denoted by~$\On$.
The small Greek letters usually stand for the ordinals.
Moreover, the following abbreviations are in common parlance:
$ \alpha <\beta\!:=\alpha\in\beta$,
$\alpha\le\beta\!:=(\alpha\in\beta)\vee (\alpha =\beta)$, and
$\alpha +1\!:=\alpha\cup\{\alpha\}$.
If $\alpha <\beta$ then $\alpha$
 {\it preceeds\/}  $\beta$ and $\beta$
 {\it succeds\/} ~$\alpha$.

Appealing to the axiom of foundation, we easily
come to a simpler definition:
an ordinal is a~transitive set totally ordered
by membership:
 $$
 \Ord(X)\leftrightarrow\Tr(X)\wedge (\forall\,u\in X)
 (\forall\,v\in X)
 (u\in v\ \vee\ u=v\ \vee\ v\in u).
 $$

\noindent
Note the useful auxiliary facts:

 {\bf(1)}
 Let $X$ and $Y$ be some classes. If $X$ is an~ordinal class,
 is transitive, and $X\ne Y$ then
 $Y\subset X$ if and only if~$Y\in X$.

 {\bf(2)}
 The intersection of two ordinal classes is an ordinal class too.

 {\bf(3)}
 If $X$ and $Y$ are ordinal classes then $ X\in Y\vee X=Y\vee Y\in X$.

 {\bf(4)} The following hold:
\begin{itemize}
\item[\bf(a)] the members of an ordinal class are ordinals;

 \item[\bf(b)] $\On$ is the only ordinal class other than any ordinal;

 \item[\bf(c)]if $\alpha$ is an ordinal then $\alpha +1$ is the
 least ordinal among the successors of~$\alpha$;

 \item[\bf(d)]the union $\bigcup X$ of a nonempty class of ordinals
 $X\subset\On$
 is an ordinal class too; if $X$ is a~set then
 $\bigcup X$ is the least upper bound for~$X$ in the ordered
 class~$\On$.
\end{itemize}

 {\bf(5)}
 The least upper bound of a~set of ordinals $x$ is usually denoted by~$\lim(x)$.
 An ordinal $\alpha$ is
 {\it limit\/} provided that
 $\alpha\ne\varnothing$ and $\lim(\alpha)=\alpha$. Equivalently,
  $\alpha$ is a~limit ordinal whenever $\alpha$ admits no presentation
  in the form $\alpha =\beta +1$ with any~$\beta\in\On$.
The symbol $K_{\rm{II}}$ stands for the class of all limit ordinals.
The ordinals, residing beyond~$K_{\rm{II}}$,
comprise the class of {\it nonlimit\/} ordinals
 $
 K_{\rm I}\!:=\On\setminus K_{\rm{II}}= \{\alpha\in\On \mid(\exists\,\beta\in\On)
 (\alpha =\beta +1)\}.
 $
Denote by~$\omega$ the least limit ordinal.
It gives pleasure to recall that the  members of~$\omega$
are {\it positive integers}.

{\bf(6)} Since
$\Ord(x)$ is a~restricted formula for~$x\in \mathbbm V$, it follows from
4.1\,(4) that
$$\alpha\in\On\leftrightarrow{\mathbbm V} ^{(B)}\models\Ord(\alpha
 ^{\scriptscriptstyle\wedge}).
 $$

 {\bf(7)}
  ${\mathbbm V}^{(B)}\models\Ord(x)$ holds for~$x\in{\mathbbm V}^{(B)}$
 if and only if there are an ordinal $\beta\in\On$ and a partition
 of unity $(b_\alpha)_{\alpha\in\beta}\subset B$ such that $x=\mix_{\alpha\in\beta}(b_\alpha\alpha
^{\scriptscriptstyle\wedge})$.

%\vskip-5pt
 \subsection{}{\bf Cardinals.}
 The concept of cardinal stems from the palaeolithic
 trick of counting by comparison with an agreed assemblage.
 This trick is witnessed  by many archeological bullas with tokens.

 Two sets are {\it equipollent\/} or have the same {\it cardinality\/}
 provided that there is a~bijection from one of them onto
 the other. An ordinal not equipollent to any of its predecessors
 is a~{\it cardinal}.
 In other words,  the cardinals make the scale
 of standards for comparing cardinalities and serve as
 referents  of the measure of quantity.
 Each positive integer is a~cardinal.
 Clearly, $\omega$ is the least infinite cardinal.
 To each infinite ordinal $\alpha$
 there is a~unique least cardinal $\mathscr H(\alpha)$
 greater than~$\alpha$. These resulting cardinals are usually called
 {\it alephs},  and their symbolizing  customarily involves the traditional
 and somewhat excessive  notations:
 $$
 \gathered[1.5\jot]
 \aleph_0:=\omega_0:=\omega;
 \\
 %\quad
 \aleph_{\alpha+1}:=\omega_{\alpha+1}:=\mathscr H(\omega_\alpha)
 \quad (\alpha\in \On);\\
 \aleph_{\beta}:=\omega_{\beta}:=\lim\{\omega_{\alpha}\mid \alpha\in\beta\}
 \quad (\beta\in K_{\rm II}).
 \endgathered
 $$

 {\bf(1)}
 Each set~$x$ is equipollent with the unique
  cardinal~$|x|$ called the  {\it cardinal number\/} or
 {\it cardinality\/} of~$x$.
 A set~$x$
 is {\it countable} in case $|x|=\omega _0\!:=\omega$ and
 {\it at most countable}, in case $|x|\le\omega _0$.

 {\bf(2)}
 The standard names of ordinals and cardinals are
 {\it standard ordinals\/} and {\it standard cardinals}.
If $\alpha\in\On$ then
 ${\mathbbm V}^{(B)}\models  \Ord(\alpha^{\scriptscriptstyle\wedge})$
 and so there is
 a unique aleph
  $\aleph_{\alpha^{\scriptscriptstyle\wedge}}$ inside  ${\mathbbm V}^{(B)}$.
As mentioned in~7.2\,(7), an arbitrary ordinal inside ${\mathbbm V}^{(B)}$
 is a~mixing of standard  ordinals.  Similarly, each Boolean valued
cardinal is a~mixing of some family of standard  cardinals.

{\bf(3)}
 $
 [\![(\omega_\alpha)^{\scriptscriptstyle\wedge}\leq
 \aleph_{\alpha^{\scriptscriptstyle\wedge}}]\!]=\mathbb 1
 $
for every cardinal  $\alpha$.

{\bf(4)} A~Boolean algebra~$B$ enjoys the {\it countable chain condition\/}
provided that every disjoint subset of~$B$ is at most countable.
In such an algebra we have
 ${\mathbbm V}^{(B)}\models
 (\omega_\alpha)^{\scriptscriptstyle\wedge}=
 \aleph_{\alpha^{\scriptscriptstyle\wedge}}
 $
for all $\alpha\in\On$.

 \subsection{}{\bf The continuum hypothesis.}
 Given an ordinal~$\alpha$, denote by~$2^{\omega _\alpha}$
 the cardinality of the powerset
  $\mathscr P (\omega _\alpha)$; i.e.,
  $
  2^{\omega _\alpha}\!:=|\mathscr P (\omega _\alpha)|.
  $
 This notation is justified since
 $2^x$ and $\mathscr P(x)$ are equipollent for all~$x$, where ${2}^x$
 is the class of all mappings from~$x$ to~$2$.
 Cantor discovered and demonstrated that $|x|<|2^x|$ for every set~$x$.
In particular, $\omega _\alpha < 2^{\omega _\alpha}$ for each ordinal~$\alpha$.
By definition, we see  that $\omega _{\alpha +1}\le 2^{\omega _\alpha}$.

 The question of whether or not there are some intermediate
 cardinalities between $\omega _{\alpha +1}$ and ${2}^{\omega _\alpha}$;
 i.e., whether or not  $\omega _{\alpha +1}= 2^{\omega _\alpha}$, is the content of
 the
 {\it generalized problem of the continuum}.
 In case $\alpha =0$, this is the classical {\it problem of the continuum}.

 The {\it continuum hypothesis\/} is in common parlance
 the equality $\omega _1 = 2^\omega $.
 The continuum hypothesis enables us to
 well order every segment of the straight line so that
 every subsegment with respect to the new order will be at most
 countable.
 The absence of intermediate cardinalities, expressed as
$(\forall \alpha\in\On)\,\omega_{\alpha +1}=2^{\omega _\alpha}$,
is the {\it generalized continuum hypothesis}.

 \subsection{}{\bf Algebras of forcing.}~
 We turn now  to describing some special class of
 Boolean algebras $B$ that opens up many amazing possibilities
 such as ``gluing together'' two standard infinite
 cardinals as well as various shifts along the alephic scale
 inside an~appropriate  Boolean valued universe $\mathbbm V^{(B)}$.

 {\bf(1)}~Consider an ordered set
 $P\!:=(P,\leq)$. Assume that $P$ has the least element~$\mathbb 0$.
 This involves no loss in generality since in case
 $P$ lacks the bottom, we can always supply the latter by
 replacing~$P$ with~$P\cup\{\mathbb 0\}$.
 We define the binary relation~$\perp$ on~$P$ by letting
 $$
 p\perp q\,\leftrightarrow\,(\forall\, r\in P)(r\leq p\wedge r\leq q
 \rightarrow r=\mathbb 0).
 $$
 Consider the {\it polar\/}
 $
 A^\perp:=\pi_{\perp}(A):=\{q\mid  (\forall p\in A)\ q\perp p\}
 $
 and use the abbreviation
 $[p]\!:=\{p\}^{\perp\perp}$. The relation $\perp$ is symmetric.
 Moreover, if $p\perp p$ then $p=\mathbb 0$.  In particular, the least
 $\perp$-band $P^\perp$
 coincides with~$\{\mathbb 0\}$. It is easy to check that
 $p\mapsto [p]$ is a monotonic mapping and $\perp$ is
 a~{\it disjointness\/} on~$P$. By the general properties of
 disjointness, the inclusion ordered set
 $\mathfrak K_{\perp}(P):=\{A^\perp\mid A\subset P\}$ of all
 $\perp$-bands in~$P$ is a~complete Boolean algebra. Recall that a~subset
 $P$ of~$B$ is {\it dense\/} provided that
 to each nonzero $b\in B$  there is a nonzero $p\in P$ satisfying
 $p\leq b$.

 {\bf(2)} For an ordered set
 $P$ the following are equivalent:
\begin{itemize}
 \item[\bf(a)]
 if $p,q\in P$ and $\mathbb 0\ne q\nleq$ then
 there is a nonzero $p^{\kern1pt\prime}\in P$ satisfying
 $p^{\kern1pt\prime}\leq q$ and $p\perp p^{\kern1pt\prime}$;

 \item[\bf(b)]
 $[p]=[\mathbb 0,p]$ for all $p\in P$;

 \item[\bf(c)]
  $p\mapsto[p]$ is one-to-one;

 \item[\bf(d)]
 $p\mapsto[p]$ is an order isomorphism of~$P$ onto a~dense subset
 of the complete Boolean algebra~$\mathfrak K_{\perp}(P)$.
\end{itemize}

An ordered set $P$ is
 {\it refined\/} provided that $P$ enjoys one
 (and, hence, all) of the conditions
 (a)--(d) of~7.5\,(2). In other words, the refined ordered sets are, up to
 isomorphism,   the dense subsets of complete Boolean algebras.
 The complete Boolean algebra $\mathfrak K_{\perp}(P)$ is usually referred to
 as the {\it Boolean completion\/} of~$P$.
 We proceed now to the examples of refined ordered sets of use in the sequel.

 {\bf (3)}~Take two nonempty sets $x$ and~$y$.
Denote by $C(x,y)$ the set of all
 functions with domain a~finite subset of~$x$ and range in~$y$.
 In other words,
 $$
 C(x,y)\!:=\{f\mid \Fnc(f)\ \wedge\ \dom(f)\in\mathscr P_{\Fin}(x)\
 \wedge\ \im(f)\subset y\}.
 $$
 Equip $C(x,y)$ with the following order
 $
 g\leq f\leftrightarrow g\supset f.
 $
 If $g\nleq f$ then
 $\dom(f)\subset\dom(g)$ and $f\ne g|_{\dom(f)}$ or, otherwise,
 $\dom(f)$ does not lie in $\dom(g)$.

 In the first case, put
 $f^\prime\!:=g$; and, in the second case, define $f^\prime$
 by the rules $\dom(f^\prime)\!:=\dom(f)\cup\dom(g)$,
 $f^\prime|_{\dom(g)}=g$, and $f^\prime|_z\ne f|_z$, where
 $z\!:=\dom(f)\setminus\dom(g)$. In both cases it is easy that
 $f^\prime\leq g$ and $f\perp f^\prime$ in $C(x,y)\cup\{\mathbb 0\}$. Hence,
  $C(x,y)\cup\{\mathbb 0\}$ is a~refined ordered set.
 We denote the  Boolean completion of~$C(x,y)$ by~$B(x,y)$.

 {\bf(4)}~Let $\varkappa$ be an~infinite cardinal. Assume that
  $x$ and~$y$ are the same as in~(3), while $|y|\geq 2$. Denote
 by $C_\varkappa(x,y)$ the set of all functions
 whose domain lies in $x$ and has cardinality strictly
 less than~$\varkappa$. In other words,
 $$
  C_\varkappa (x,y)\!:=\{f \mid  \Fnc(f)\ \wedge\ \dom(f)\in\mathscr P(x)\
 \wedge\ |\dom(f)|<\varkappa\wedge\im(f)\subset y\}.
 $$
 Order $C _\varkappa(x,y)$ in the same manner as in~(3).
 Arguing as above, we see that
 $C_\varkappa (x,y)\cup\{\mathbb 0\}$ is a~refined ordered set.
 We denote the  Boolean completion of~$C _\varkappa(x,y)$  by
 $B_\varkappa(x,y)$. Clearly,
 $C(x,y)=C_\omega (x,y)$ and $B(x,y)=B_\omega (x,y)$.

\subsection{}{\bf Shifting alephs.}
The collection of Boolean algebras in~7.5 enables us to achieve
a cardinal shift in an~appropriate Boolean valued model, so
demonstrating Cohen's beautiful result.
In his talk of 2001 on the discovery of forcing,
Cohen mentioned that ``using the language of Boolean
algebras brings our technique of forcing
to standard usages''~\cite[p.~1096]{Cohen-2002}.

{\bf(1)}
 Let $\lambda$ be an arbitrary infinite cardinal. Consider
 the complete Boolean algebra ~$B\!:=B(\omega,\lambda)$. Then
 $|\lambda^{\scriptscriptstyle\wedge}|$ is the countable cardinal inside
 $\mathbbm V^{(B)}$; i.e.,
  $
  \mathbbm  V^{(B)}\models|\lambda^{\scriptscriptstyle\wedge}|=\aleph_0.
  $

{\bf(2)}
Given two infinite  cardinals
 $\varkappa$ and $\lambda$, there is a complete
 Boolean algebra~$B$ such that
 $$
 \mathbbm V^{(B)}\models\!|\varkappa^{\scriptscriptstyle\wedge}|
 =|\lambda^{\scriptscriptstyle\wedge}|.
 $$

{\bf(3)}
$B(x,2)$ is a Boolean algebra with the countable chain condition.

{\bf(4)}
 Let $x$ be a~nonempty set and $|x|=\omega_\alpha$.
 Then
$$
 \omega_\alpha\leq|B(x,2)|\leq(\omega_\alpha)^{\omega_0}.
 $$

{\bf(5)}
If  $(\omega_\alpha)^{\omega_0}=\omega_\alpha$
and $B\!:=B(\omega\times\omega_\alpha,2)$ then
$
\mathbbm V^{(B)}\models\!2^{\aleph_0}=
\aleph_{\alpha^{\scriptscriptstyle\wedge}}.
$

{\bf(6)} If $\ZF$ is consistent then
$\ZFC$ remains consistent together with the axiom $2^{\,\omega_0}=\omega_2$.

Indeed, G\"odel demonstrated~\cite{Goedel} that
$\ZF$ remains consistent on assuming the axiom of choice and the generalized continuum hypothesis.
Hence, we may suppose that
$$ (\omega_2)^{\omega_0}=\big(2^{\,\omega _1}\big)^{\omega_0}=
 2^{\,\omega_1\cdot\omega_0}=2^{\,\omega_1}=\omega_2.
$$
By~(5) there is a~complete Boolean
 algebra~$B$ such that $\mathbbm V^{(B)}\models\!2^{\aleph_0}=
 \aleph_{2^{\scriptscriptstyle\wedge}}$. Since
  $\mathbbm V^{(B)}\models\!2^{\scriptscriptstyle\wedge}=2$,
it is obvious that
 $\mathbbm V^{(B)}\models\!\aleph_{2^{\scriptscriptstyle\wedge}}=\aleph_2$
 and so $\mathbbm V^{(B)}\models\!2^{\aleph_0}=\aleph_2$. Moreover,
 $\mathbbm V^{(B)}\models\ZFC$ by the transfer principle 3.2\,(1).

{\bf(7)}
Continuity of all homomorphisms of the classical
Banach algebra $C([0,1])$
is independent of the rest of the axioms of~$\ZFC$ but depends on
the cardinality of~$2^{\omega_0}$.
On this matter, see~\cite[p.~19]{Truth}.

{\bf(8)}
Let $A$ stand for the set of functions from $\mathbbm R$ to
countable subsets of~$\mathbbm R$.
If $f\in A$ while $x$ and $y$ are some random reals
then it is reasonable to believe that
$y\in f(x)$ with probability~0; i.e., $x\notin f(y)$
with probability~1. Similarly, $y\notin f(x)$ with probability~1.
Using these instances of the common sense,
Freiling  proposed~\cite{Freiling}  the  axiom $\AX$ which reads
$$
(\forall f\in A)(\exists x)(\exists y)\ x\notin f(y)\wedge y\notin f(x).
$$
It turns out that  $\AX$ amounts to the negation of
the continuum hypothesis,  $\CH$, in~$\ZFC$.
This is another argument in favor
of the Cohen statement: ``I think the consensus will be that
$\CH$ is false''~\cite[p.~1099]{Cohen-2002}.

The above approach to the problem of the continuum
proceeds mainly  along the lines of~\cite{Bell}.
As regards the omitted details and references, see~\cite[Ch.~9]{IBA}.
The state of the art in forcing is presented in~\cite{Shelah}.

\section{\bf Logic and Freedom}

   Mathematics is  the most ancient of sciences. However, in the beginning was the
   word. We must remember that the olden ``logos'' lives in logics and
   logistics rather than grammar. The order of mind and the order of store
   are the precious gifts of our ancestors.

   The intellectual field  resides beyond the grips of the law of
   diminishing returns. The more we know, the huger become the frontiers
   with the unbeknown, the oftener we meet the mysterious. The twentieth
   century enriched our geometrical views with the concepts of space-time
   and fractality.  Each instance of knowledge is an event, a point in the
   Minkowski 4-space. The realm of our knowledge comprises a clearly bounded
   set of these instances. The frontiers of science produce the boundary between the
   known and the unknown which is undoubtedly fractal and we have no grounds
   to assume it rectifiable or measurable.
   It is worth noting in parentheses that rather smooth are the routes to the
   frontiers of science
   which are charted by teachers, professors, and all other kinds
   of educationalists.
      Pedagogics  dislikes saltations and sharp changes of the prevailing
   paradigm. Possibly, these topological obstructions reflect some objective
   difficulties in modernizing education.

    The proofs are uncountable of the fractality of the
    boundary between the known and the unbeknown. Among them we
    see such negative trends as the unleashed growth of pseudoscience, mysticism,
    and other forms of obscurantism which creep into all lacunas of the
    unbeknown.   As revelations of fractality appear the most unexpected,
    beautiful, and stunning interrelations between
    seemingly distant areas and directions of
    science.

   The revolutionary changes in mathematics  in the vicinity of
   the turn
   of the twentieth century are connected not only with
   the new calculus of the infinite which was created by Cantor.
   Of similar import was the rise and development
   of mathematical logic which applied  rigor and analysis to
   the very process of mathematical demonstration.  The decidable and
   the undecidable, the provable and the improvable, the consistent
   and the inconsistent have entered the research lexicon of the
   perfect mathematician.  Mathematics became a reflexive
   science  that is engaged  not only in search of truths but also
   in study of its own methods for attaining these truths.

   Aristotle's logic, the paradoxes of Zeno, the razor of William of
   Occam, the donkey of Buridan,  the  Lebnizian   Calculemus, and
    Boolean  algebras are the outstanding achievements of  mankind
    which cast light on the road to the new stages of logical
    studies. Frege immortalized his name by inventing the calculus of
    predicates which underlies the modern logic.

   The twentieth century is marked with deep
   penetration of the ideas of mathematical logic into
   many sections of science and technology. Logic is a tool that
   not only organizes and orders our ways of thinking but also
   liberates us from dogmatism in choosing the objects and methods
   of  research.  Logic of today is  a major instrument and
   instituion of mathematical freedom.
   Boolean valued analysis serves as a brilliant confirmation of this thesis.

   Returning to Takeuti's original definition of Boolean valued analysis,
   we must acknowledge its extraordinary breadth.
   The Boolean valued model resting on  the dilemma of ``verum'' and ``falsum''
   is employed implicitly by the overwhelming majority of mathematicians.
   Our routine talks and discussions on seminars hardly deserve the qualification
   of articles of prose. By analogy, it seems pretentious to claim that
   Euler, Cauchy, and Abel exercised  Boolean valued analysis.

    Boolean valued analysis is a special mathematical technique
    based on  validating truth by means of a nontrivial
    Boolean algebra.
    From a category-theoretic viewpoint,
    Boolean valued analysis  is the theory of Boolean toposes.
    From a topological viewpoint,  it is the theory of
    continuous polyverses over Stone spaces.

   Mach taught us the economy of thought. It seems reasonable to
   apply his principle and to shorten the bulky term
  ``Boolean valued analysis.''
   Mathematization of the laws of thought originated with
   Boole \cite{Boole}  and deserves the lapidary title ``Boolean analysis.''
\endgroup

%\newpage

\bibliographystyle{plain}

\begin{thebibliography}{99}

\bibitem{Takeuti}
Takeuti G. (1978)
{\it Two Applications of Logic to Mathematics.}
Tokio--Princeton: Iwanami Publ. \& Princeton University Press.


\bibitem{Scott}
Scott D. (1969)
``Boolean Models and Nonstandard Analysis,''
In:  {\it Applications of Model Theory to Algebra, Analysis, 
and Probability} (Ed.: Luxemburg W.~A.~J.). 
New York etc.: Holt, Rinehart, and Winston, 87--92.


\bibitem{Cohen}
Cohen P. J. (1966)
{\it Set Theory and the Continuum Hypothesis.} New York etc.:  Benjamin.

\bibitem{Problems}
Hilbert D. (1902)
``Mathematical problems.
Lecture delivered before the International Congress of Mathematicians
at Paris in 1900,''
{\it Bull. Amer. Math. Soc.}, {\bf8}, 437--479.

 \bibitem{Johnstone}
 Johnstone~P.~T. (2002)
 {\it Sketches of an Elephant. A~Topos Theory Compendium.}
  Oxford: Clarendon Press (Oxford Logic Guides;~438).


\bibitem{Hilbert}
Hilbert D. (1967)
``On the infinite,'' In: {\it From Frege to G\"odel 1879--1931:
A Source Book in the History of Science}. Cambridge: Harvard University Press,
367--392.



\bibitem{Hofstadter}
Hofstadter~D.~R. (1999)
{\it G\"odel, Escher, Bach: an~Eternal Golden Braid
(20th Anniversary Edition)}. New York: Basic Books.

 \bibitem{Gordon}
Gordon~E.~I. (1977)
 ``Real numbers in Boolean-valued models
of set theory, and $K$-spaces,''
{\it Soviet Math. Doklady}, {\bf18}, 1481--1484.



\bibitem{IBA}
Kusraev A.~G. and Kutateladze S.~S. (2005)
{\it Introduction to Boolean Valued Analysis}.
Moscow:  Nauka Publishers (in Russian).

\bibitem{Luzin}
Luzin~N.~N. (1928)
``Real function theory: the state of the art,''
In: {\it Proceedings of the All-Russia Congress of Mathematicians (Moscow,  April 27--May 4, 1927)},
 Moscow and Leningrad: Glavnauka, 11--32.


\bibitem{Goedel}
G\"odel K.  (1940)
{\it The Consistency of the Axiom of Choice and
of the Generalized Continuum Hypothesis}. Princeton: Princeton Univ. Press.

\bibitem{Yandell}
Yandell B.~H. (2002)
{\it The Honors Class. Hilbert's Problems and Their Solvers.}
Natick: A.~K.~Peters, Ltd.

\bibitem{Kanamori}
Kanamori A. (1996)
``The mathematical developments of set theory
from Cantor to Cohen,'' Bull. Symbolic Logic, {\bf1}:1, 1--70.

\bibitem{Cohen-2002}
Cohen P. (2002)
 ``The discovery of forcing,'' Rocky Mountain J.
 Math., {\bf32}:4, 1071--1100.

\bibitem{Manin}
Manin Yu.~I. (2004)
``Georg Cantor and his heritage,''
Proceedings of the Steklov Institute, {\bf246}, 208--216.

\bibitem{Truth}
Dales~H.~G. and Oliveri~G. (Eds.) (1998)
{\it Truth in Mathematics.}
Oxford: Clarendon Press.

\bibitem{Freiling}
Freiling~Ch. (1986)
``Axioms of symmetry: throwing darts at the real line,''
J. Symbolic Logic, {\bf51}, 190--200.



\bibitem{Bell}
Bell~J.~L. (2005)
{\it Set Theory: Boolean-Valued Models and Independence Proofs.}
Oxford: Clarendon Press (Oxford Logic Guides;~47).

\bibitem{Shelah}
Shelah~S. (1998)
{\it Proper and Improper Forcing.}
Berlin: Springer-Verlag.

\bibitem{Boole}
Boole~G. (1997)
{\it Selected Manuscripts on Logic and Its Philosophy}.
 Basel: Birkh\"auser-Verlag (Science Networks. Historical Studies;~{20}).




\end{thebibliography}

\end{document}